\pdfoutput=1
\documentclass[11pt]{article}

\usepackage[utf8]{inputenc}
\usepackage[T1]{fontenc}
\usepackage[margin=1in]{geometry}
\usepackage{amsmath}
\usepackage{amsthm}
\usepackage{amsfonts}
\usepackage{amssymb}
\usepackage{graphicx}
\usepackage{booktabs}
\usepackage{multirow}
\usepackage{diagbox}
\usepackage{array}
\usepackage{rotating}
\usepackage{makecell}
\usepackage{xcolor}
\usepackage{colortbl}
\usepackage{bm}
\usepackage{siunitx}
\usepackage{float}
\usepackage{algorithm}
\usepackage{algpseudocode}
\usepackage{algorithmicx}
\usepackage{adjustbox}
\usepackage[numbers,sort&compress]{natbib}
\usepackage{hyperref}
\usepackage[capitalize,nameinlink]{cleveref}

\colorlet{blue}{black}
\colorlet{red}{black}

\hypersetup{
    colorlinks=true,
    linkcolor=black,
    citecolor=black,
    urlcolor=black,
}

\theoremstyle{plain}
\newtheorem{theorem}{Theorem}
\newtheorem{lemma}{Lemma}
\theoremstyle{remark}
\newtheorem{remark}{Remark}

\floatname{algorithm}{Algorithm}

\makeatletter
\providecommand{\theHALG@line}{}
\renewcommand{\theHALG@line}{\thealgorithm.\arabic{ALG@line}}
\makeatother

\makeatletter
\let\original@includegraphics\includegraphics
\renewcommand{\includegraphics}[2][]{%
  \IfFileExists{#2}{%
    \original@includegraphics[#1]{#2}%
  }{%
    \IfFileExists{#2.pdf}{%
      \original@includegraphics[#1]{#2}%
    }{%
      \IfFileExists{#2.png}{%
        \original@includegraphics[#1]{#2}%
      }{%
        \IfFileExists{#2.jpg}{%
          \original@includegraphics[#1]{#2}%
        }{%
          \fbox{\begin{minipage}[c][0.28\textheight][c]{0.88\linewidth}
          \centering Missing figure file: \texttt{\detokenize{#2}}\\
          Include this file in the arXiv source package to display the figure.
          \end{minipage}}%
        }%
      }%
    }%
  }%
}
\makeatother

\newcommand{\keywords}[1]{\par\medskip\noindent\textbf{Keywords:} #1\par\bigskip}

\begin{document}

\title{An efficient and memory free algorithm for subdiffusion equation using incremental singular value decomposition}

\author{Jichun Li\thanks{Department of Mathematical Science, University of Nevada, Las Vegas, Las Vegas, Nevada. Email: \texttt{Jichun.li@unlv.edu}.}
\and Yangpeng Zhang\thanks{Department of Mathematical Science, University of Nevada, Las Vegas, Las Vegas, Nevada. Corresponding author. Email: \texttt{yangpeng.zhang@unlv.edu}.}
\and Yangwen Zhang\thanks{Department of Mathematics, University of Louisiana at Lafayette, Lafayette, Louisiana. Email: \texttt{yangwen.zhang@louisiana.edu}.}}

\date{\today}
\maketitle

\begin{abstract}
In this paper, we address the well-known challenge in the numerical solution of time-fractional partial differential equations (TFPDEs), namely, that the dependence on all previous time levels leads to storage requirements that grow linearly with the number of time steps. To overcome this difficulty, we develop an efficient algorithm based on incremental singular value decomposition (ISVD), which avoids the excessive memory demands associated with storing the full solution history. A rigorous error analysis is established, and numerical experiments are presented to validate the theoretical results. Comparisons with the direct method and a representative fast evaluation method show that the proposed ISVD approach dramatically reduces memory usage relative to the direct method and remains competitive with the fast method over the tested parameter regimes.
\end{abstract}

\keywords{Time-fractional partial differential equation; Caputo derivative; Incremental singular value decomposition; Finite element method}

\section{Introduction}\label{sec:intro}
Time-fractional partial differential equations (TFPDEs) generalize classical evolution equations by incorporating fractional-order Caputo temporal derivatives. These models capture non-local memory effects and anomalous diffusion phenomena absent in integer-order counterparts, with transformative applications in contaminant transport through heterogeneous media, protein diffusion in crowded cellular environments, and charge carrier dynamics in disordered semiconductors  \cite{Benson2000,METZLER20001,Scher1975}. Time-fractional PDEs  introduce historical dependence that imposes severe computational burdens: storage requirements grow linearly with time steps  \cite{Kilbas2006,Garrappa2018}. 

The limitation is particularly acute for high-dimensional problems and long-time simulations, necessitating advanced numerical strategies. Seminal works  \cite{Jin2013} established optimal $L^2$-error estimates for linear FEM on convex domains. For temporal discretization, the $L1$ scheme remains predominant due to its balance of accuracy and stability  \cite{Li_SISC2011,Zeng2013, SakamotoYamamoto2011,Mustapha2018, Le2016}. An analysis of robust Crank-Nicolson scheme for the sub-diffusion equation is derived by  \cite{JinLiZhou2016CN}. Li et al.  \cite{LiYangZhou2023splitFEM} present a splitting method for the sub-diffusion equation, applying high-order FEM on the smooth part and tailored discretization for the nonsmooth part, achieving high-order accuracy. Lin and Xu  \cite{LIN20071533} proved its $O\left(\tau^{2-\alpha}\right)$ convergence for smooth solutions. Ren et al.  \cite{Ren2018fast,Ren2017superconvergence} presented two fully discrete schemes for diffusion equations and diffusion-wave equations, respectively. The fully discrete system requires storing all historical solutions, creating heavy memory demands for practical-scale simulations. 

In this paper, we consider the following sub-diffusion equation:
\begin{equation}\label{Fractional_PDEs}
\left\{
\begin{aligned}
\partial_{t}^{\alpha} u(\mathbf{x},t) - \Delta u(\mathbf{x},t) &= f(\mathbf{x},t) && (\mathbf{x},t) \in \Omega \times (0,T], \\
u(\mathbf{x},t) &= 0 && (\mathbf{x},t) \in \partial \Omega \times (0,T], \\
u(\mathbf{x},0) &= u_0(\mathbf{x}) && \mathbf{x} \in \Omega,
\end{aligned}
\right.
\end{equation}
where $\Omega$ is a bounded convex polytopal domain in $\mathbb{R}^{d}~ (d=2,3)$ with a boundary $\partial \Omega$ and $u_0(\mathbf{x})$ is a given function defined on the domain $\Omega$ and $T>0$ is the final time. Here $\partial_{t}^{\alpha} u~(0<\alpha<1)$ denotes the left-sided Caputo fractional derivative of order $\alpha$ with respect to $t$ and it is defined by (see  \cite{Kilbas2006}, p.70)
$$
\textcolor{red}{\partial_{t}^{\alpha} u(\mathbf{x},t)=\frac{1}{\Gamma(1-\alpha)} \int_{0}^{t}(t-s)^{-\alpha} \frac{\mathrm{\partial}}{\mathrm{\partial} s} u(\mathbf{x},s) \mathrm{d} s},
$$
where $\Gamma(\cdot)$ is the standard Gamma function defined by $\Gamma(x)=\int_{0}^{\infty} s^{x-1} \mathrm{e}^{-s} \mathrm{~d} s$ for $\Re(x)>0 .$

First we describe the spatial discretization based on the standard Galerkin finite element method (FEM). Let $\mathcal{T}_{h}$ be a shape regular and quasiuniform triangulation of domain $\Omega$ into $d$-simplexes, known as finite elements and denoted by $K$. Then over the triangulation $\mathcal{T}_{h}$ we define a continuous piecewise linear finite element space $V_{h}$ by
$$
V_{h}=\left\{v_{h} \in H_{0}^{1}(\Omega):\left.v_{h}\right|_{K} \text { is a linear function, } \forall K \in \mathcal{T}_{h}\right\},
$$
where $H_0^1(\Omega)$ is a standard Sobolev space:
\[
H_0^1(\Omega) = \left\{u\in L^2(\Omega): Du\in \left(L^2(\Omega)\right)^d, u=0 ~\text{on} ~\partial\Omega\right\}.
\]
Then the semidiscrete Galerkin scheme for problem \eqref{Fractional_PDEs} reads: find $u_{h}(t) \in V_{h}$ such that
\begin{equation}\label{semiFEM}
\left(\partial_{t}^{\alpha} u_{h}, \chi\right)+\left(\nabla u_{h}, \nabla \chi\right)=(f, \chi) \quad \forall \chi \in V_{h}, t>0.
\end{equation}
For the time discretization, we divide the interval $[0, T]$ into $n$ equally spaced subintervals with a time step size $\tau=T / n$, and $t_{i}=i \tau, i=0, \ldots, n$. Then the $L1$ scheme  \cite{LIN20071533,Li_SISC2011} approximates the Caputo fractional derivative $\partial_{t}^{\alpha} u\left(\mathbf{x}, t_{i}\right)$ by
$$
\begin{aligned}
\partial_{t}^{\alpha} u\left(\mathbf{x}, t_{i}\right) &=\frac{1}{\Gamma(1-\alpha)} \sum_{j=0}^{i-1} \int_{t_{j}}^{t_{j+1}} \frac{\partial u(\mathbf{x}, s)}{\partial s}\left(t_{i}-s\right)^{-\alpha} \mathrm{d} s \\
& \approx \frac{1}{\Gamma(1-\alpha)} \sum_{j=0}^{i-1} \frac{u\left(\mathbf{x}, t_{j+1}\right)-u\left(\mathbf{x}, t_{j}\right)}{\tau} \int_{t_{j}}^{t_{j+1}}\left(t_{i}-s\right)^{-\alpha} \mathrm{d} s \\
&=\sum_{j=0}^{i-1} \beta_{j} \frac{u\left(\mathbf{x}, t_{i-j}\right)-u\left(\mathbf{x}, t_{i-j-1}\right)}{\tau^{\alpha} \Gamma(2-\alpha)}=: \bar{\partial}_{\tau}^{\alpha} u\left(t_{i}\right), 
\end{aligned}
$$
where the weights $\left\{\beta_{j}\right\}$ are given by
$$
\beta_{j}=(j+1)^{1-\alpha}-j^{1-\alpha}, \quad j=0,1, \ldots, i-1.
$$
Then the fully discrete scheme reads: given $u_{h}^{0}= P u_0 \in V_{h}$ with $L^2$-projection operator $P: L^2(\Omega)\longrightarrow V_h$. Find $u_{h}^{i} \in V_{h}$ for $i=1,2, \ldots, n$ such that
\begin{align}\label{FEM}
	\begin{split}
&\beta_{0}\left(u_{h}^{i}, v_h\right)+c_{\alpha} \tau^\alpha\left(\nabla u_h^{i}, \nabla v_h\right)\\
&=\beta_{i-1}\left(u_h^{0}, v_h\right)+\sum_{j=1}^{i-1}\left(\beta_{j-1}-\beta_{j}\right)\left(u_h^{i-j}, v_h\right)+c_{\alpha} \tau^{\alpha}\left(f\left(t_{i}\right), v_h\right)\textcolor{red},
\end{split}
\end{align}
where $c_{\alpha}=\Gamma(2-\alpha)$. The computational challenge of the fully discrete scheme is obvious: To compute the numerical solution $u_{h}^{i}$, the solutions $ \{u_{h}^{j} \}_{j=0}^{i-1}$ at all preceding time instances are required, as a result of the nonlocality of the Caputo fractional derivative $\partial_{t}^{\alpha} u$. The direct method requires $\mathcal{O}(mn)$ storage and $\mathcal{O}(mn^2)$ computational work, where $m$ denotes the number of spatial degrees of freedom and $n$ denotes the number of time steps. Consequently, while the storage requirement grows linearly with $n$, the overall computational cost grows quadratically with $n$. This becomes a significant bottleneck, particularly for high-dimensional problems and multi-query applications.  Here we provide \Cref{fractionalPDEs} which is standard $L1$ scheme as the direct method to update the numerical solution. For the sake of readability, we consistently use bold formatting for matrices and vectors throughout our work. For example, $\mathbf{u}_0$ in \Cref{fractionalPDEs} is a column vector that contains all coefficients of the numerical solution $u_h^0$ for some basis functions $\{\phi_s\}_{s=1}^m$ at first time level after some proper interpolation $P$ for initial condition $u_0(\mathbf{x})$ that  will be specified later. $\mathbf{M}$ and $\mathbf{A}$ are standard mass and stiffness matrices respectively.

\begin{algorithm}[H]
\caption{ (Standard finite element and $L1$ scheme for solving the time-fractional partial differential equation)}
\label{fractionalPDEs}
{\bf Input:} $c_\alpha$, $\alpha$, $\tau$, $\mathbf{M}\in \mathbb R^{m\times m}$, $\mathbf{A}\in \mathbb R^{m\times m}$, $\mathbf{u}_0\in\mathbb{R}^{m}$
\begin{algorithmic}[1]
	\State  Set $\widetilde{\mathbf{u}}_0  = \mathbf{M}\mathbf{u}_0$; $\widetilde{\mathbf{A}} = \mathbf{M} + c_{\alpha} \tau^{\alpha} \mathbf{A}$; $\bm \beta = \texttt{zeros}(n,1)$; $\widetilde{\mathbf{U}} = \texttt{zeros}(m,n)$
	\For {$i=1$ to $n$} 
	\State $\bm \beta (i) = i^{1-\alpha} - (i-1)^{1-\alpha}$;
	\State Get the load vector at $t_i$, denoted by $\mathbf{b}_i$;
	\State $\mathbf{b}_i \gets\bm  \beta(i)\widetilde{\mathbf{u}}_0 + c_{\alpha} \tau^{\alpha}  \mathbf{b}_i + \sum_{j=1}^{i-1} \widetilde{\mathbf{U}}(:,i-j) (\bm\beta(j-1) - \bm \beta(j))$;
	\State Solve $\widetilde{\mathbf{A}}\mathbf{u}_i = \mathbf{b}_i$ to get $\mathbf{u}_i$;
	\State $\widetilde{\mathbf{U}}(:,i) \gets \mathbf{Mu}_i$;
	\EndFor
\end{algorithmic}
{\bf Output:} $\mathbf{u}_n$
\end{algorithm}

The computational solution of time-fractional differential equations is inherently challenging due to their non-local operators, which require storing and integrating the entire solution history. To mitigate the associated prohibitive memory and computational costs, fast evaluation method using sum-of-exponentials (SOE) approximations  \cite{Jiang2015, Zhang_fast_evaluation2017,Schadle2006} has emerged as a standard approach. This method decomposes the temporal convolution integral into a local part and a history part. The local part is typically discretized using standard $L1$ schemes, while the history part leverages an SOE approximation of the singular kernel $t^{-1-\alpha}$. The primary advantage of SOE methods lies in their ability to approximate the history contribution using only $N_{exp}$ exponential terms, where $N_{exp}$ scales logarithmically $(O(\log N_T)$ or $O(\log^{2} N_T))$ with the total number of time steps $N_T$, for a given accuracy tolerance $\varepsilon$. This drastically reduces storage requirements to $O(N_S N_{exp})$ and computational work to $O(N_S N_T N_{exp})$, compared to $O(N_S N_T)$ storage and $O(N_S N_T^{2})$ work for direct methods ($N_S$ being the number of spatial grid points).

{\color{blue} Beyond SOE-based fast evaluation, a substantial literature has developed fast-history algorithms for time-fractional problems from several complementary viewpoints, including fast time stepping for fractional integrals  \cite{Li2010fasttime}, interval-clustering methods for evolution equations with memory  \cite{McLean2012interval}, fast and oblivious convolution quadrature  \cite{Schadle2006}, second-order accurate fast schemes for subdiffusion  \cite{Zeng2015second}, high-dimensional fast difference schemes  \cite{Zeng2016fasthighdim}, and efficient high-order convolution-quadrature algorithms  \cite{BanjaiLopez2019}. These methods accelerate the treatment of the nonlocal history term through kernel approximation, clustering, or fast convolution, and they provide important points of comparison for the present work.}

In this paper, we introduce  an efficient methodology based on incremental singular value decomposition (ISVD) to resolve this issue when substantial history storage is required. The foundational work introduced by Brand \cite{Brand2006} gives the ISVD-based low-rank updates. Rather than the fast evaluation algorithm by using approximation for the kernel $t^{-1-\alpha}$, we keep the temporal discretization form in $L1$ type. Our approach comprises two key parts: 
\begin{itemize}
\item Low-rank Compression: Solution states generated during time stepping are dynamically compressed using low-rank approximation techniques. This exploits the inherent structure and redundancy often present in the evolving solution field.
\item Efficient History Storage: The compressed representations are stored efficiently using SVD, significantly reducing the memory footprint associated with retaining the full solution history.
\end{itemize}

This framework achieves substantial reductions in both storage complexity and computational workload while rigorously maintaining solution accuracy. Crucially, a major advantage of our ISVD-based approach is its generality. The compression mechanism is agnostic to the specific temporal discretization scheme employed for the TFPDE. This flexibility renders the methodology broadly applicable to diverse classes of history-dependent TFPDEs, irrespective of their chosen time-stepping strategy.

{\color{blue} In the numerical section, we compare the proposed approach with the direct method and with a representative SOE-based fast-history solver that uses the same $L1$/finite element framework, so that the comparison is carried out under identical spatial and temporal discretizations. Accordingly, our performance claims are stated in a restricted form: the proposed ISVD method achieves large memory savings relative to the direct method and is competitive with the fast method in the tested regimes.}

The paper is organized as follows. In Section 2 we describe the ISVD algorithm in detail and provide the algorithms. In Section 3, we apply ISVD to solve the sub-diffusion equation and present the computational cost. In Section 4, we present an error analysis for the $L1$ scheme for the sub-diffusion equation. {\color{blue} In Section 5, we present several fractional diffusion PDE experiments to compare ISVD with the direct method and a representative fast evaluation algorithm. The numerical results show that ISVD preserves the expected convergence behavior, dramatically reduces memory usage relative to the direct method, and remains competitive with the fast method in the tested parameter regimes.} Throughout the paper, we denote by $C$ a generic constant, which may vary at different occurrences but is always independent of the mesh size $h$, the solution $u$, and the initial data $u_0$.

\section{The Incremental Singular Value Decomposition Method}\label{sec3}

This section presents fundamental definitions and the incremental Singular Value Decomposition (SVD) algorithm. For $\mathbf{u} \in \mathbb{R}^n$ and $r \leq n$, $\mathbf{u}(1:r)$ denotes the first $r$ components of $\mathbf{u}$. Similarly, for $\mathbf{U} \in \mathbb{R}^{m \times n}$, $\mathbf{U}(p:q, r:s)$ refers to the submatrix from rows $p$ to $q$ and columns $r$ to $s$. We assume $\mathbf{U}$ is low-rank: $\texttt{rank}(\mathbf{U}) \ll \min\{m, n\}$.

We adapt Brand's incremental SVD algorithm from  \cite{Zhang2022isvd} for updating the SVD when columns are appended. The process consists of four steps:

\subsection*{Step 1: Initialization}
Given a nonzero initial column $\mathbf{u}_1 \in \mathbb{R}^m$, initialize the SVD as:
\begin{align*}
  \Sigma = \|\mathbf{u}_1\|_2, \quad \mathbf{Q} = \mathbf{u}_1 \Sigma^{-1}, \quad R = 1,
\end{align*}
where $\|\cdot\|_2$ is standard Euclidean norm $\|\mathbf{u}\|_2=\left(u_1^2+u_2^2+\cdots u_m^2\right)^{1/2}$. The following Algorithm \ref{algorithm01} formalizes this step.

\begin{algorithm}[H]
  \caption{Initialize Incremental SVD (ISVD)}
  \label{algorithm01}
  \textbf{Input:} $\mathbf{u}_1 \in \mathbb{R}^m$ 
  \begin{algorithmic}[1]
    \State $\Sigma \gets (\mathbf{u}_1^\top \mathbf{u}_1)^{1/2}$ 
    \State $\mathbf{Q} \gets \mathbf{u}_1 \Sigma^{-1}$  
    \State $R \gets 1$
  \end{algorithmic}
  \textbf{Output:} $\mathbf{Q}, \Sigma, R$
\end{algorithm}

\subsection*{Step 2: $p$-Truncation for Sequential Updates}
Assume the rank-$k$ truncated SVD of the first $\ell$ columns $\mathbf{U}_\ell$ is:
\begin{align}
  \mathbf{U}_\ell \approx \mathbf{Q}\mathbf{\Sigma} \mathbf{R}^\top, \quad
  \mathbf{Q}^\top \mathbf{Q} = \mathbf{I}_{k}, \quad
  \mathbf{R}^\top \mathbf{R} = \mathbf{I}_{k}, \quad
  \mathbf{\Sigma} = \text{diag}(\sigma_1,\dots,\sigma_k). \label{eq201}
\end{align}
Suppose the next $s$ incoming columns satisfy
 \textcolor{blue}{ \begin{align}
\|\mathbf{u}_i - \mathbf{Q}\mathbf{Q}^\top \mathbf{u}_i\|_2 &< \mathtt{tol1}, \quad i = \ell+1,\dots,\ell+s, \label{lessthantol} \\
  \|\mathbf{u}_i - \mathbf{Q}\mathbf{Q}^\top \mathbf{u}_i\|_2 &\geq \mathtt{tol1}, \quad i = \ell+s+1. \label{largerthantol}
\end{align}}
For columns satisfying \eqref{lessthantol}, update $\mathbf{U}_{\ell+s} = [\mathbf{U}_\ell \mid \mathbf{u}_{\ell+1} \mid \cdots \mid \mathbf{u}_{\ell+s}]$:
\begin{align*}
  \mathbf{U}_{\ell+s} &\approx \mathbf{Q} \mathbf{Y} \begin{bmatrix} \mathbf{R}&\mathbf{0} \\\mathbf{0} & \mathbf{I}_s \end{bmatrix}^\top, \quad
  \mathbf{Y} = [\mathbf{\Sigma} \mid \mathbf{Q}^\top \mathbf{u}_{\ell+1} \mid \cdots \mid \mathbf{Q}^\top \mathbf{u}_{\ell+s}],
\end{align*}
where top-right $\mathbf{0}\in\mathbb{R}^{k\times s}$ and bottom-left $\mathbf{0}\in\mathbb{R}^{s\times k}$. Compute the SVD $\mathbf{Y} = \mathbf{Q}_\mathbf{Y} \mathbf{\Sigma}_\mathbf{Y} \mathbf{R}_\mathbf{Y}^\top$ and partition $\mathbf{R}_\mathbf{Y} = \begin{bmatrix} \mathbf{R}_\mathbf{Y}^{(1)} \\ \mathbf{R}_\mathbf{Y}^{(2)} \end{bmatrix}$. Update:
\begin{align*}
  \mathbf{Q} &\gets \mathbf{Q} \mathbf{Q}_\mathbf{Y}, \quad
  \mathbf{\Sigma} \gets \mathbf{\Sigma}_\mathbf{Y}, \quad
  \mathbf{R} \gets \begin{bmatrix} \mathbf{R} \mathbf{R}_\mathbf{Y}^{(1)} \\ \mathbf{R}_\mathbf{Y}^{(2)} \end{bmatrix}.
\end{align*}
The projection matrix $\mathbf{W} = [\mathbf{Q}^\top \mathbf{u}_{\ell+1} \mid \cdots \mid \mathbf{Q}^\top \mathbf{u}_{\ell+s}]$ is stored incrementally.

\subsection*{Step 3: Full Update for Significant Residuals}
For $\mathbf{u}_{\ell+s+1}$ satisfying \eqref{largerthantol}, compute the residual:
 \textcolor{blue}{ \begin{align}
  \mathbf{e} = \mathbf{u}_{\ell+s+1} - \mathbf{Q}\mathbf{Q}^\top \mathbf{u}_{\ell+s+1}, \quad p = \|\mathbf{e}\|_2 \geq \mathtt{tol1}, \quad \widetilde{\mathbf{e}} = \mathbf{e}/p. \label{residual}
\end{align}}
Update $\mathbf{U}_{\ell+s+1} = [\mathbf{U}_{\ell+s} \mid \mathbf{u}_{\ell+s+1}]$:
\begin{align*}
  \mathbf{U}_{\ell+s+1} &\approx [\mathbf{Q} \mid \widetilde{\mathbf{e}}] \bar{\mathbf{Y}} \begin{bmatrix} \mathbf{R}^\top & \mathbf{0} \\ \mathbf{0}^\top & 1 \end{bmatrix}, \quad
  \bar{\mathbf{Y}} = \begin{bmatrix} \mathbf{\Sigma} & \mathbf{Q}^\top \mathbf{u}_{\ell+s+1} \\ \mathbf{0}^\top & p \end{bmatrix}.
\end{align*}
Compute the SVD $\bar{\mathbf{Y}} = \bar{\mathbf{Q}} \bar{\mathbf{\Sigma}} \bar{\mathbf{R}}^\top$ and update:
\begin{align*}
  \mathbf{Q} &\gets [\mathbf{Q} \mid \widetilde{\mathbf{e}}] \bar{\mathbf{Q}}, \quad
  \mathbf{\Sigma} \gets \bar{\mathbf{\Sigma}}, \quad
  \mathbf{R} \gets \begin{bmatrix} \mathbf{R} & \mathbf{0} \\ \mathbf{0}^\top & 1 \end{bmatrix} \bar{\mathbf{R}}.
\end{align*}
The matrix dimensions increase in this step.

\begin{remark}
  {Orthogonality of $\mathbf{e}$ to the column space of $\mathbf{Q}$ may be lost numerically \cite{Giraud2005}. To maintain orthogonality, two iterative refinement steps are recommended when $|\mathbf{e}^\top \mathbf{Q}(:,1)| > \mathtt{tol2}$ \cite{Zhang2022isvd}.}
\end{remark}

\subsection*{Step 4: Singular Value Truncation}
\begin{lemma}
Assume that $\mathbf{\Sigma} = diag (\sigma_1,\sigma_2,\cdots,\sigma_k)$ with $\sigma_1 \ge \sigma_2 \ge\cdots \ge\sigma_k$, and $\bar{\mathbf{\Sigma}} = \text{diag}(\mu_1, \dots, \mu_{k+1})$ with $\mu_1 \geq \cdots \geq \mu_{k+1}$, we have
\begin{align}
  \mu_{k+1} &\leq p, \label{eq302} \\
  \mu_{k+1} &\leq \sigma_k \leq \mu_k \leq \sigma_{k-1} \leq \cdots \leq \sigma_1 \leq \mu_1. \label{eq303}
\end{align}
\end{lemma}
\begin{proof}
See Lemma 3, in  \cite{Zhang2022isvd}.
\end{proof}
Inequality \eqref{eq302} implies the last singular value can be small, necessitating truncation. By \eqref{eq303}, only the last singular value requires checking:
\begin{itemize}
   \item  \textcolor{blue}{ If $\bar{\mathbf{\Sigma}}(k+1,k+1) \geq \mathtt{tol3}$}, retain full SVD.
  \item \textcolor{blue}{  If $\bar{\mathbf{\Sigma}}(k+1,k+1) < \mathtt{tol3}$}, truncate:
  \begin{align*}
    \mathbf{Q} &\gets [\mathbf{Q} \mid \widetilde{\mathbf{e}}] \bar{\mathbf{Q}}(:,1:k), \\
    \mathbf{\Sigma} &\gets \bar{\mathbf{\Sigma}}(1:k,1:k), \\
    \mathbf{R} &\gets \begin{bmatrix} \mathbf{R} & \mathbf{0} \\ \mathbf{0}^\top & 1 \end{bmatrix} \bar{\mathbf{R}}(:,1:k).
  \end{align*}
\end{itemize}

\begin{lemma}\label{error_SingularValueTruncation}
  Let $\mathbf{Q}\mathbf{\Sigma}\mathbf{R}^\top$ be the SVD of $\mathbf{A} \in \mathbb{R}^{m \times n}$ with singular values $\sigma_1 \geq \cdots \geq \sigma_r > 0$. Define $\mathbf{B} = \mathbf{Q}(:,1:r-1)\mathbf{\Sigma}(1:r-1,1:r-1)\mathbf{R}(:,1:r-1)^\top$. Then:
  \begin{align*}
    \max_{1 \leq i \leq n} \|\mathbf{a}_i - \mathbf{b}_i\|_2 \leq \sigma_r,
  \end{align*}
  where $\mathbf{a}_i$ and $\mathbf{b}_i$ are the $i$-th columns of $\mathbf{A}$ and $\mathbf{B}$, and $\|\cdot\|_2$ denotes the Euclidean norm.
\end{lemma}
The proof follows directly from the properties of SVD. Algorithm \ref{alg202} summarizes the update procedure.
\section{Incremental SVD method  for the time fractional PDE}\label{ISVD_IDeq}
	This section focuses on the application of the incremental SVD algorithm to the sub-diffusion equation \eqref{Fractional_PDEs}.  Our approach is to simultaneously solve the TFPDEs and incrementally update the SVD of the solution. By doing so, we store the solutions at all time steps in the four matrices of the SVD. As a result, we are able to address the issue of data storage in solving the TFPDEs \eqref{Fractional_PDEs}. 
	
	Due to the errors that may arise during the data compression process and the potential alterations caused by singular value truncation to previous storage, it becomes necessary for us to modify the traditional scheme \eqref{FEM}. Below, we provide a brief discussion of our approach.
\begin{algorithm}[H]
  \caption{Update Incremental SVD}
  \label{alg202}
  \textbf{Input:} $\mathbf{Q} \in \mathbb{R}^{m \times k}$, $\mathbf{\Sigma} \in \mathbb{R}^{k \times k}$, $\mathbf{R} \in \mathbb{R}^{\ell \times k}$,  \textcolor{blue}{ $\mathtt{tol1}$, $\mathtt{tol2}$, $\mathtt{tol3}$}, $\mathbf{W}$, $\mathbf{Q}_0$, $q$, $\mathbf{u}_{\ell+1}$
  \begin{algorithmic}[1]
    \State $\mathbf{d} \gets \mathbf{Q}^\top \mathbf{u}_{\ell+1}$
    \State $\mathbf{e} \gets \mathbf{u}_{\ell+1} - \mathbf{Q}\mathbf{d}$; $p \gets \|\mathbf{e}\|_2$
 	 \If{ \textcolor{blue}{ $p \geq \mathtt{tol1}$}}
      \If{$q > 0$} 
        \State $\mathbf{Y} \gets [\mathbf{\Sigma} \mid \text{cell2mat}(\mathbf{W})]$
        \State $[\mathbf{Q}_\mathbf{Y}, \mathbf{\Sigma}_\mathbf{Y}, \mathbf{R}_\mathbf{Y}] \gets \texttt{svd}(\mathbf{Y}, \text{'econ'})$
        \State $\mathbf{Q}_0 \gets \mathbf{Q}_0 \mathbf{Q}_\mathbf{Y}$; $\mathbf{\Sigma} \gets \mathbf{\Sigma}_\mathbf{Y}$
        \State $\mathbf{R}_1 \gets \mathbf{R}_\mathbf{Y}(1:k,:)$; $\mathbf{R}_2 \gets \mathbf{R}_\mathbf{Y}(k+1:\text{end},:)$
        \State $\mathbf{R} \gets \begin{bmatrix} \mathbf{R}\mathbf{R}_1 \\ \mathbf{R}_2 \end{bmatrix}$; $\mathbf{d} \gets \mathbf{Q}_\mathbf{Y}^\top \mathbf{d}$
      \EndIf
      \State $\bar{\mathbf{Y}} \gets \begin{bmatrix} \mathbf{\Sigma} & \mathbf{d} \\ \mathbf{0}^\top & p \end{bmatrix}$
      \State $[\bar{\mathbf{Q}}, \bar{\mathbf{\Sigma}}, \bar{\mathbf{R}}] \gets \texttt{svd}(\bar{\mathbf{Y}})$; $\widetilde{\mathbf{e}} \gets \mathbf{e}/p$
      \If{ \textcolor{blue}{ $|\widetilde{\mathbf{e}}^\top \mathbf{Q}(:,1)| > \mathtt{tol2}$}} \Comment{Orthogonality refinement}
        \State $\mathbf{e} \gets \mathbf{e} - \mathbf{Q}(\mathbf{Q}^\top \mathbf{e})$; $p_1 \gets \|\mathbf{e}\|_2$; $\widetilde{\mathbf{e}} \gets \mathbf{e}/p_1$
      \EndIf
      \State $\mathbf{Q}_0 \gets \begin{bmatrix} \mathbf{Q}_0 & \mathbf{0} \\ \mathbf{0}^\top & 1 \end{bmatrix} \bar{\mathbf{Q}}$
      \If{ \textcolor{blue}{ $\bar{\mathbf{\Sigma}}(k+1,k+1) \geq \mathtt{tol3}$}}
        \State $\mathbf{Q} \gets [\mathbf{Q} \mid \widetilde{\mathbf{e}}] \mathbf{Q}_0$; $\mathbf{\Sigma} \gets \bar{\mathbf{\Sigma}}$; $\mathbf{R} \gets \begin{bmatrix} \mathbf{R} & \mathbf{0} \\ \mathbf{0}^\top & 1 \end{bmatrix} \bar{\mathbf{R}}$
        \State $\mathbf{Q}_0 \gets \mathbf{I}_{k+1}$
      \Else
        \State $\mathbf{Q} \gets [\mathbf{Q} \mid \widetilde{\mathbf{e}}] \mathbf{Q}_0(:,1:k)$; $\mathbf{\Sigma} \gets \bar{\mathbf{\Sigma}}(1:k,1:k)$
        \State $\mathbf{R} \gets \begin{bmatrix} \mathbf{R} & \mathbf{0} \\ \mathbf{0}^\top & 1 \end{bmatrix} \bar{\mathbf{R}}(:,1:k)$; $\mathbf{Q}_0 \gets \mathbf{I}_{k}$
      \EndIf
      \State $\mathbf{W} \gets [\,]$; $q \gets 0$
    \Else
      \State $q \gets q + 1$; $\mathbf{W}\{q\} \gets \mathbf{d}$
    \EndIf
  \end{algorithmic}
  \textbf{Output:} $\mathbf{Q}, \mathbf{\Sigma}, \mathbf{R}, \mathbf{W}, \mathbf{Q}_0, q$
\end{algorithm}

\begin{remark}
	\textcolor{blue}{In Algorithm \ref{alg202}, the index $k$ denotes the current number of retained singular values, i.e., the current retained rank of the compressed representation. It is not a time-step index. The three tolerances play different roles: $\mathtt{tol1}$ is the threshold used to decide whether the new snapshot is already well represented by the current column space, $\mathtt{tol2}$ is the numerical threshold used to monitor loss of orthogonality and trigger reorthogonalization, and $\mathtt{tol3}$ is the singular-value truncation threshold. In practice, $\mathtt{tol1}$ and $\mathtt{tol3}$ should be chosen below the target discretization error, while $\mathtt{tol2}$ is typically taken near machine precision. In the benchmark problems reported in this paper, we use $\mathtt{tol1}=\mathtt{tol2}=\mathtt{tol3}=10^{-12}$ for most cases. }
\end{remark}
\begin{remark}
	\textcolor{blue}{If $\mathbf{W}$ is non-empty, the output of Algorithm \ref{alg202} is not yet the exact SVD of $\mathbf{U}_{\ell+1}$. In this case, $\mathbf{W}$ stores the coefficient vectors of snapshots that are already well represented by the current basis and are therefore postponed to a later SVD update rather than expanded immediately. Thus, the information is not discarded. The compressed representation $\mathbf{Q}[\mathbf{\Sigma} \mathbf{R}^{\top} \mid \mathbf{W}]$ is still the quantity used in the subsequent time-stepping procedure, and the postponed data in $\mathbf{W}$ are absorbed in the next update when a later snapshot satisfies $p \ge \mathtt{tol1}$.}
\end{remark}

\begin{enumerate}
\item \textbf{Initialization and Time Step 1:}
Begin with the initial condition $u_h^0$. Compute the numerical solution at $t_1$ using traditional methods, denoted as $\widehat{u}_h^1$. Compress the solution history ${u_h^0, \widehat{u}_h^1}$ to obtain the first compressed dataset ${\widetilde{u}_h^{1,0}, \widetilde{u}_h^{1,1}}$, where:
\begin{itemize}
\item $\widetilde{u}_h^{1,0}$ represents the compressed initial condition,
\item $\widetilde{u}_h^{1,1}$ represents the compressed solution at $t_1$.
\end{itemize}
\item \textbf{Time Step 2:}
Using the compressed data $\{\widetilde{u}_h^{1,0}, \widetilde{u}_h^{1,1}\}$, compute the numerical solution at $t_2$, denoted as $\widehat{u}_h^2$. Compress the extended solution history $\{\widetilde{u}_h^{1,0}, \widetilde{u}_h^{1,1}, \widehat{u}_h^2\}$ to obtain $\{\widetilde{u}_h^{2,j}\}_{j=0}^2$, where:
\begin{itemize}
    \item $\widetilde{u}_h^{2,0}, \widetilde{u}_h^{2,1}$ are updated compressed representations of earlier solutions,
    \item $\widetilde{u}_h^{2,2}$ is the new compressed solution at $t_2$.
\end{itemize}

\item \textbf{General Time Step} $\mathbf{i}$:
Given compressed data $\{\widetilde{u}_h^{i,j}\}_{j=0}^i$ representing the solution history up to $t_i$:
\begin{enumerate}
    \item Compute the numerical solution at $t_{i+1}$, denoted as $\widehat{u}_h^{i+1}$,
    \item Compress the augmented solution set $\{\widetilde{u}_h^{i,j}\}_{j=0}^i \cup \widehat{u}_h^{i+1}$ to generate updated compressed data $\{\widetilde{u}_h^{i+1,j}\}_{j=0}^{i+1}$.
\end{enumerate}
Here, $\widetilde{u}_h^{i+1,j}$ denotes the compressed representation of the solution at $t_j$ using basis functions available at $t_{i+1}$.

\item \textbf{Final Time Step:}
Continue the process until reaching the terminal time $t_N$. The final output is $\widehat{u}_h^N$, with the complete compressed history $\{\widetilde{u}_h^{N,j}\}_{j=0}^N$ available for analysis.
\end{enumerate}
In summary, we apply our proposed approach by incrementally compressing data at each time step to compute the numerical solutions throughout the process.
	
Based on the preceding discussion, we can present our formulation below, where we seek $\widehat{u}_h^{i} \in V_h$ that satisfies the following equation:
	\begin{align}\label{ISVD_eq1}
		\begin{split}
		&\beta_{0}\left(\widehat{u}_{h}^{i}, v_h\right)+c_{\alpha} \tau^\alpha\left(\nabla \widehat {u}_h^{i}, \nabla v_h\right)\\
		&=\beta_{i-1}\left(\widetilde {u}_h^{i-1,0}, v_h\right)+\sum_{j=1}^{i-1}\left(\beta_{j-1}-\beta_{j}\right)\left(\widetilde{u}_h^{i-1,i-j}, v_h\right)+c_{\alpha} \tau^{\alpha}\left(f\left(t_{i}\right), v_h\right).
		\end{split}
	\end{align}
	
Subsequently, we express equation \eqref{ISVD_eq1} into matrix form to highlight the benefits of utilizing the incremental SVD for solving the TFPDEs more distinctly.  To do this, let $\widehat {\mathbf{u}}_{i}$ and $\widetilde {\mathbf{u}}_{i-1,j}$ denote the coefficient of the basis functions for $\widehat {{u}}_h^{i}$ and $\widetilde {{u}}_h^{i-1,j}$ and $\mathbf{b}_i$ the coefficient of the basis function after interpolation for source term $f$ at time level $i$, respectively. Then we seek a solution $\widehat {\mathbf{u}}_{i}\in \mathbb R^m$ that satisfies the following equation:
	\begin{align}\label{ISVD_eq2}
		(\beta_0\mathbf{M}+c_\alpha \tau^\alpha \mathbf{A})\widehat {\mathbf{u}}_{i}=&\mathbf{M} \sum_{j=1}^{i-1} (\beta_{j-1} - \beta_j) \widetilde {\mathbf{u}}_{i-1,i-j}\\\notag
		&+\beta_{i-1}\mathbf{M}\widetilde {\mathbf{u}}_{i-1, 0}+c_\alpha \tau^\alpha \mathbf{b}_i.
	\end{align}
	Here, $\{\widetilde {\mathbf{u}}_{i-1,j}\}_{j=0}^{i-1}$  represents the data that has been compressed from \\$\{\widetilde{\mathbf{u}}_{i-2, 0}, \ldots, \widetilde{\mathbf{u}}_{i-2, i-2}, \widehat {\mathbf{u}}_{i-1}\}$ using the incremental SVD algorithm. We assume that $\mathbf{Q}_{i-1}$, $\mathbf{\Sigma}_{i-1}$, $\mathbf{R}_{i-1}$, and $\mathbf{W}_{i-1}$ are the matrices associated with this compression process. In other words,
	
\begin{align}\label{ISVD_eq3}
[\widetilde{\mathbf{u}}_{i-2,0}\mid&\widetilde{\mathbf{u}}_{i-2,1}\mid \cdots\mid \widetilde{\mathbf{u}}_{i-2,i-2}\mid \widehat{\mathbf{u}}_{i-1}]\\\notag
&\xrightarrow[\text{}]{\text{Compress}}
\mathbf{Q}_{i-1} [\mathbf{\Sigma}_{i-1}\mathbf{R}_{i-1}\mid \mathbf{W}_{i-1}]^\top
= [\widetilde{\mathbf{u}}_{i-1,0}\mid \cdots\mid \widetilde{\mathbf{u}}_{i-1,i-1}] .
\end{align}
	
	Let $[\mathbf{\Sigma}_{i-1}\mathbf{R}_{i-1}\mid \mathbf{W}_{i-1}]^\top$ be denoted as $\mathbf{X}_{i-1}$.  Accordingly, equation \eqref{ISVD_eq2} can be rewritten as follows:
	\begin{align}\label{ISVD_eq4}
		(\beta_0\mathbf{M}+c_\alpha \tau^\alpha \mathbf{A})\widehat {\mathbf{u}}_{i}=&\mathbf{M} \mathbf{Q}_{i-1}\sum_{j=1}^{i-1} (\beta_{j-1} - \beta_j) \mathbf{X}_{i-1}(:,i-j)\\\notag
		&~+\beta_{i-1}\mathbf{M}\widetilde{\mathbf{ u}}_{i-1, 0}+c_\alpha\tau^\alpha \mathbf{b}_i.
	\end{align}
	
\begin{figure}[tbh]
		\centerline{
			\hbox{\includegraphics[width=\textwidth]{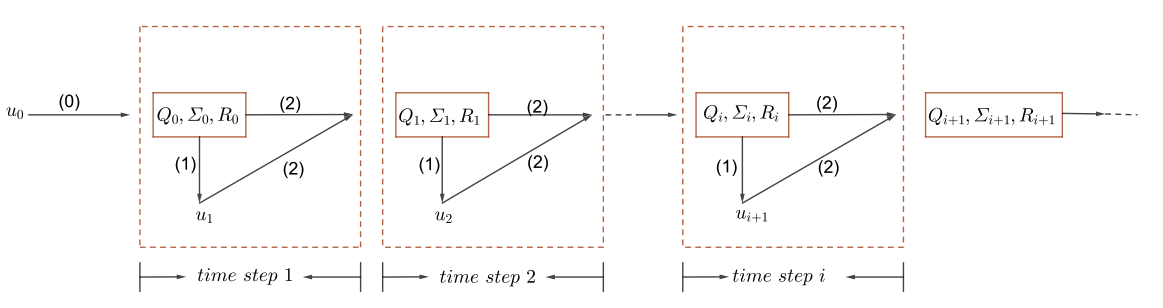}}}
		\caption{The process of using the incremental SVD to solve the time fractional PDE.}
		\label{ISVD_figure}
		\centering
	\end{figure}

	Once $\widehat{\mathbf{u}}_{i}$ is obtained, we update the SVD of $[\widetilde {\mathbf{u}}_{i-1,0}\mid \cdots\mid \widetilde {\mathbf{u}}_{i-1,i-1}\mid \widehat {\mathbf{u}}_{i}]$ using $\mathbf{Q}_{i-1}$, $\mathbf{\Sigma}_{i-1}$, $\mathbf{R}_{i-1}$, $\mathbf{W}_{i-1}$, and $\widehat {\mathbf{u}}_{i}$ based on the principles of the incremental SVD method. This update process is illustrated  in  \Cref{ISVD_figure}.
	
	
    Throughout the remainder of this section, we will examine the memory and computational cost pertaining to the history term in our proposed approach. Our data storage involves four matrices: $\mathbf{Q}_{i-1}$, $\mathbf{\Sigma}_{i-1}$, $\mathbf{R}_{i-1}$, and $\mathbf{W}_{i-1}$, resulting in a memory cost of $\mathcal O((m+n)r)$, where $r$ represents the rank of the solution data. By taking into account our assumption that $r\ll \min\{m, n\}$, we can compare this memory cost with that of the traditional approach presented in \eqref{fractionalPDEs}, which shows the reduction achieved by the proposed approach.
	
	Moving on to the computational cost, which also encompasses the cost of the incremental SVD, it can be summarized as follows:
	\begin{align*}
		\mathcal O(mnr) + \sum_{i=1}^n \sum_{j=0}^{i} \mathcal O(r) = \mathcal O(mnr + rn^2).
	\end{align*}
	
	Here, once again, $r$ represents the rank of the solution data. Based on our assumption that $r\ll \min\{m, n\}$, we can compare the computational cost in \eqref{FEM} to that of the traditional approach, revealing that our approach experiences only linear growth, rather than quadratic growth as observed in the traditional approach.
	
	\begin{algorithm}[H]
\caption{Solve Sub-diffusion equation with Incremental SVD}
\label{fpdeandisvd}
\textbf{Input:} 
$c_\alpha$, $\alpha$, $\tau$, $n$, 
$\mathbf{M} \in \mathbb{R}^{m \times m}$, $\mathbf{A} \in \mathbb{R}^{m \times m}$, 
$\mathbf{u}_0 \in \mathbb{R}^m$,  \textcolor{blue}{ $\mathtt{tol1}$,$\mathtt{tol2}$,$\mathtt{tol3}$}

\begin{algorithmic}[1]
\State Initialize: 
    $\mathbf{V} \gets [\,]$, $\mathbf{Q}_0 \gets 1$, $q \gets 0$, 
    $\widetilde{\mathbf{u}}_0 \gets \mathbf{M}\mathbf{u}_0$, 
    $\widetilde{\mathbf{A}} \gets \mathbf{M} + c_{\alpha} \tau^{\alpha} \mathbf{A}$, 
    $\bm{\beta} \gets \texttt{zeros}(n, 1)$
    
\State $[\mathbf{Q}, \bm{\Sigma}, \mathbf{R}] \gets \texttt{InitializeISVD}(\mathbf{u}_0)$  
\For{$i = 1$ \textbf{to} $n$}  
    \State $\bm{\beta}(i) \gets i^{1-\alpha} - (i-1)^{1-\alpha}$  \Comment{Compute fractional coefficient}
    \State Obtain load vector $\mathbf{b}_i$ at time $t_i$
    
    \If{$q = 0$}  
        \State $\texttt{MQ} \gets \mathbf{M}\mathbf{Q}$  \Comment{Precompute for efficiency}
        \State $\mathbf{C} \gets \bm{\Sigma} \mathbf{R}^\top$  
        \State $\mathbf{b}_i \gets \bm{\beta}(i)\widetilde{\mathbf{u}}_0 + c_{\alpha}\tau^{\alpha}\mathbf{b}_i + 
                 \texttt{MQ} \left( \sum_{j=1}^{i-1} \mathbf{C}(:, i-j) (\bm{\beta}(j-1)-\bm \beta(j)) \right)$
        \State Solve $\widetilde{\mathbf{A}} \mathbf{u}_i = \mathbf{b}_i$ for $\mathbf{u}_i$
    \Else  
        \State $\mathbf{D} \gets [\mathbf{C} \mid \texttt{cell2mat}(\mathbf{V})]$  
        \State $\mathbf{b}_i \gets \bm{\beta}(i)\widetilde{\mathbf{u}}_0 + c_{\alpha}\tau^{\alpha}\mathbf{b}_i + 
                 \texttt{MQ} \left( \sum_{j=1}^{i-1} \mathbf{D}(:, i-j) (\bm{\beta}(j-1)-\bm \beta(j))  \right)$
        \State Solve $\widetilde{\mathbf{A}} \mathbf{u}_i = \mathbf{b}_i$ for $\mathbf{u}_i$
    \EndIf
    
    \State Update ISVD: 
        $[q, \mathbf{V}, \mathbf{Q}_0, \mathbf{Q}, \bm{\Sigma}, \mathbf{R}] \gets$ \\
        ~~~~~~~~~~~\hspace*{0.5cm}\texttt{UpdateISVD}$(q, \mathbf{V}, \mathbf{Q}_0, \mathbf{Q}, \bm{\Sigma}, \mathbf{R}, \mathbf{u}_i,  \textcolor{blue}{ \mathtt{tol1},\mathtt{tol2},\mathtt{tol3}})$
    
    \If{$q > 0$}  
        \State $\mathbf{C} \gets \bm{\Sigma} \mathbf{R}^\top$
    \EndIf
\EndFor
\end{algorithmic}
\textbf{Output:} Solution vector $\mathbf{u}_n$
\end{algorithm}

\begin{remark}
\color{blue}	Recent works have employed incremental SVD techniques to compress simulation data and to update singular triplets or reduced bases; see, for example, Fareed et al.~ \cite{MR3775096,Areed2020}, Li et al.~ \cite{li2024incremental}, and Costanzo et al.~ \cite{MR4982783}. The role of ISVD in the present paper is different. In our approach, after each time step, the newly computed finite element solution is simultaneously compressed, the full history snapshot is released from memory, and the compressed representation is used directly to evaluate the memory term at subsequent time levels in the original PDE solver. Therefore, the compressed data are not used solely for storage reduction;  instead, they are immediately fed back into the numerical time-stepping procedure. Consequently, singular-value truncation alters the reconstructed historical states that enter later solves, which makes the resulting formulation and error analysis  different from those arising in standard ISVD compression settings. To the best of our knowledge, such a mechanism that simultaneously compresses, releases, and reuses history data has not been reported in the existing literature.
\end{remark}

\section{Error estimate}\label{Error_estimate}
	In this section, we derive the error between the solution of the scheme \eqref{ISVD_eq1} and the exact solution that satisfies the equation \eqref{Fractional_PDEs}.
Now we state the main result of our paper.
\begin{theorem}\label{main_res}
Let $u$ and $\widehat {u}_h^n$ denote the solutions of \eqref{Fractional_PDEs} and \eqref{ISVD_eq1}, respectively. Throughout the entire process of the incremental SVD algorithm, the tolerance accuracy \textup{\texttt{tol1}}, \textup{\texttt{tol2}} and \textup{\texttt{tol3}} are applied to be the same as \textup{\texttt{tol}} to both $p$-truncation and singular value truncation.  Given $f$ and $u_0(x)$ sufficiently smooth, let $u$ and $u_h^n$ be the solutions of problems \eqref{Fractional_PDEs} and \eqref{FEM}, respectively.  Assume that $u\in H^4(\Omega), u_t\in H^2(\Omega), u_{tt}\in L^2(\Omega)$. Then there holds
\begin{align*}
		\left\|u\left(t_n\right)-\widehat{u}_h^n\right\|_{L^2(\Omega)} \leq C\left(h^2+\tau^{2-\alpha}+n \textup{\texttt{tol}}\right),
\end{align*}
where $C$ is a positive constant independent of $h$.
\end{theorem}

{\color{red}
	\begin{remark}
	\Cref{main_res} is stated for the uniform partition $t_i=i\tau$. Nevertheless, the same ISVD history-compression strategy can be combined with the $L1$ discretization on a graded temporal mesh $0=t_0<t_1<\cdots<t_n=T$. In that case, the uniform coefficients $\beta_j$ are replaced by the corresponding nonuniform $L1$ weights, and the fully discrete system becomes time-step dependent. Apart from this change in the discrete history weights, the compression update and the error  analysis remain essentially the same in structure. The main analytical modification is that the standard uniform-mesh $L1$ discretization estimate used below must be replaced by the corresponding graded-mesh $L1$ estimate; accordingly, the temporal term in the final error bound is replaced by the corresponding graded-mesh rate. If one wishes to treat the usual weak singularity near $t=0$ explicitly, then the standard graded-mesh regularity framework for the underlying $L1$ analysis should be adopted.
	\end{remark}
}

	\subsection[Proof of Theorem~\ref{main_res}]{Proof of \Cref{main_res}}
	We begin by giving an error bound between the solution of the standard finite element method given by equation \eqref{FEM} and the solution of the continuous problem \eqref{Fractional_PDEs}. Additionally, we derive an error bound between the solution of the standard finite element method \eqref{FEM} and our proposed scheme \eqref{ISVD_eq1}. By applying the triangle inequality, we obtain a straightforward error bound between the solution of the TFPDE \eqref{Fractional_PDEs} and our proposed scheme \eqref{ISVD_eq1}.
	
	\begin{lemma}\label{theorem401}
	Assume $f $ and $u_0$ are sufficiently smooth. Let $u$ and $u_h^n$ be the solutions of problems \eqref{Fractional_PDEs} and \eqref{FEM}, respectively. Assume that $u\in H^4(\Omega), u_t\in H^2(\Omega), u_{tt}\in L^2(\Omega)$. Then there holds
	\begin{align*}
	\left\|u\left(t_n\right)-u_h^n\right\|_{L^2(\Omega)} \leq &C\left(h^2 +\tau^{2- \alpha}\right).
	\end{align*}
	\end{lemma}
	
	\begin{proof}
	See Theorem 3.2, in  \cite{LIN20071533} and Theorem 3.16 in  \cite{Jin2015AN}. (Or see Theorem 3.5 in  \cite{Jin2013}, Theorem 3.7 in  \cite{Jin2014}, Theorem 3.16, Remark 3.17 in  \cite{Jin2015AN} for nonsmooth data $f$ and $u_0$).
	\end{proof}

	Now, we will proceed to establish the error estimation between the solution of the standard finite element method \eqref{FEM} and our proposed scheme \eqref{ISVD_eq1}.
	
	\begin{lemma}\label{error_uh_uhat}
		Let $u_h^i$ and $\widehat{u}_h^{i}$ be the solutions of \eqref{FEM} and \eqref{ISVD_eq1}, respectively. The following error bound is established:
		\begin{align*}
			\|u_{h}^{i}-\widehat{u}_{h}^{i} \| \le  \max_{0\le j\le i-1} \| u_h^j-\widetilde{u}_h^{i-1,j} \|.
		\end{align*}
	\end{lemma} 
	\begin{proof}
		Recall that $u_h^{i}$ and $\widehat  u_h^{i}$ satisfy the following equations
		\begin{subequations}
			\begin{align}
						&\beta_{0}\left( u_{h}^{i}, v_h\right)+c_{\alpha} \tau^\alpha\left(\nabla  u_h^{i}, \nabla v_h\right)\nonumber \\
				&=\beta_{i-1}\left( u_h^{0}, v_h\right)+\sum_{j=1}^{i-1}\left(\beta_{j-1}-\beta_{j}\right)\left( u_h^{i-j}, v_h\right)+c_{\alpha} \tau^{\alpha}\left(f\left(t_{i}\right), v_h\right),\label{fem_so}\\	
			&\beta_{0}\left(\widehat u_{h}^{i}, v_h\right)+c_{\alpha} \tau^\alpha\left(\nabla \widehat u_h^{i}, \nabla v_h\right)\nonumber\\
			&=\beta_{i-1}\left(\widetilde u_h^{i-1,0}, v_h\right)+\sum_{j=1}^{i-1}\left(\beta_{j-1}-\beta_{j}\right)\left(\widetilde u_h^{i-1,i-j}, v_h\right)+c_{\alpha} \tau^{\alpha}\left(f\left(t_{i}\right), v_h\right).\label{ISVD_so}
			\end{align}
		\end{subequations}
		
		Subtracting    \eqref{fem_so}  from $ \eqref{ISVD_so} $ and introducing two notations:   $\widehat{e}_{i}=u_{h}^{i}-\widehat{u}_{h}^{i}$, and $\widetilde{e}_{i-1,j}=u_{h}^{j}-\widetilde{u}_{h}^{i-1,j}$, $j=0,1,\ldots,i-1$, we  obtain the following equation:
		\begin{align}\label{eq616}
			&\beta_0\left(\widehat{e}_{i},v_h\right)+ c_{\alpha} \tau^\alpha(\nabla\widehat{e}_{i},\nabla v_h)\\\notag
			&=  \beta_{i-1} (\widetilde e_{i-1,0}, v_h)  + \sum_{j=1}^{i-1}(\beta_{j-1} - \beta_j)(\widetilde{e}_{i-1,i-j},v_h),\quad \forall v_h\in V_h.
		\end{align}
	
		Substituting $v_h = \widehat{e}_{i} \in V_h$ into the above equation  and  we derive the following equation:
		\begin{align*}
			&\beta_0\|\widehat{e}_{i}\|^2+ c_{\alpha} \tau^\alpha \|\nabla\widehat{e}_{i}\|^2\\
			&=  \beta_{i-1} (\widetilde e_{i-1,0}, \widehat{e}_{i})  + \sum_{j=1}^{i-1}(\beta_{j-1} - \beta_j)(\widetilde{e}_{i-1,i-j},\widehat{e}_{i})\\
			&\le \beta_{i-1}\frac{\|\widetilde e_{i-1,0}\|^2 + \|\widehat{e}_{i}\|^2}{2} + \sum_{j=1}^{i-1}(\beta_{j-1} - \beta_j) \frac{\|\widetilde{e}_{i-1,i-j}\|^2 + \|\widehat{e}_{i}\|^2}{2}\\
			&= \dfrac{\beta_{i-1}}{2}\|\widetilde e_{i-1,0}\|^2  +  \dfrac{\beta_{0}}{2}\|\widehat e_{i}\|^2  +  \sum_{j=1}^{i-1}(\beta_{j-1} - \beta_j) \frac{\|\widetilde{e}_{i-1,i-j}\|^2 }{2}\\
			&\le \dfrac{\beta_{0}}{2}\|\widehat e_{i}\|^2 + \dfrac{\beta_{0}}{2} \max_{0\le j\le i-1}\|\widetilde{e}_{i-1,j}\|^2.
		\end{align*}
	\end{proof}
	
	After obtaining the numerical solution $\widehat{u}_h^{i+1}$, following the approach of the incremental SVD algorithm (see \Cref{ISVD_figure}), the subsequent step involves compressing the data $\widehat{u}_h^{i+1}$ using the SVD matrices $Q_i$, $\Sigma_i$, and $R_i$. This compression process includes three cases: $p$-truncation, no-truncation, and singular value truncation. As discussed in \Cref{ISVD_IDeq}, $p$-truncation and no-truncation do not alter  the previous solution, while singular value truncation may potentially change the entire previous solution. Nevertheless, we can establish the following bound:
	\begin{lemma}\label{error_uhat_uhtilde}
		
		Let $\widehat{u}_h^{i+1}$ be the solution of \eqref{ISVD_eq1}, and let $\{\widetilde{u}_h^{i+1,j}\}_{j=0}^{i+1}$ represent the compressed solution corresponding to $\{ \widetilde{u}_h^{i,0}, \widetilde{u}_h^{i,1},\ldots, \widetilde{u}_h^{i,i}, \widehat{u}_h^{i+1}\}$. This compressed solution is obtained using the incremental SVD method with a tolerance accuracy of \textup{\texttt{tol}} applied to both $p$-truncation and singular value truncation. In this scenario, the following inequality holds:
		\begin{align*}
			\max\left\{\left\{\max_{0\le j \le i}		\|\widetilde  {u}_h^{i,j}-\widetilde{u}_h^{i+1, j}\|\right\},  \|\widehat {u}_h^{i+1}-\widetilde{u}_h^{i+1, i+1}\| \right\} \le Ch \textup{\texttt{tol}},
		\end{align*}
		where $C$ is a positive constant independent of $h$(see Theorem 3.5  \cite{LarsonBengzon2013}.).
	\end{lemma} 
	
	\begin{proof}
    		Assuming that $\widehat{\mathbf{u}}_{i+1}$ and $\widetilde{\mathbf{u}}_{k,\ell}$ are the coefficients of $\widehat {u}_h^{i+1}$ and $\widetilde u_h^{k,\ell}$ corresponding to the finite element basis functions $\{\phi_s\}_{s=1}^m$ as we mentioned before, respectively, i.e.,
\begin{equation}
\widehat {u}_h^{i+1}=\sum_{s=1}^m\hat{u}_{i+1}^s\phi_s,~~\widetilde{u}_h^{k,\ell}=\sum_{s=1}^m\widetilde{u}_{k,\ell}^s\phi_s.
\end{equation}		
where $\widehat{\mathbf{u}}_{i+1}=\left[\hat{u}_{i+1}^{1},\hat{u}_{i+1}^2,\cdots,\hat{u}_{i+1}^{m}\right]^{T},\widetilde{\mathbf{u}}_{k,\ell}=\left[\widetilde{u}_{k,\ell}^{1},\widetilde{u}_{k,\ell}^2,\cdots,\widetilde{u}_{k,\ell}^{m}\right]^{T}$. Let $\mathbf{Q}_i\mathbf{\Sigma}_i \mathbf{R}_i^\top$ be the truncated SVD of $[\mathbf{u}_0\mid \widehat{\mathbf{u}}_1\mid\ldots\mid \widehat{\mathbf{u}}_i]$ using the incremental SVD algorithm. Let $\{\widetilde{\mathbf{u}}_{i,j}\}_{j=0}^i$ be the compressed data at the $i$-th step. In other words,
		\begin{align*}
			[\mathbf{u}_0\mid \widehat{\mathbf{u}}_1\mid\ldots\mid \widehat{\mathbf{u}}_i]\xrightarrow[\text{}]{\text{Incremental SVD}} \mathbf{Q}_i\mathbf{\Sigma}_i \mathbf{R}_i^\top = 	[\widetilde{\mathbf{u}}_{i,0}\mid \widetilde{\mathbf{u}}_{i,1}\mid\ldots\mid \widetilde{\mathbf{u}}_{i,i}].
		\end{align*}
		
		Next, we will divide our discussions into three cases: $p$-truncation, no truncation, and singular value truncation.
		
		{\bf Case 1:} $p$ - truncation. At this moment, $p = \|\widehat{\mathbf{u}}_{i+1} - \mathbf{Q}_i\mathbf{Q}_i^\top \widehat{\mathbf{u}}_{i+1} \| \le \textup{\texttt{tol}}$, and 
		\begin{align*}
			[\mathbf{Q}_i\mathbf{\Sigma}_i \mathbf{R}_i^\top \mid \mathbf{Q}_i\mathbf{Q}_i^\top  \widehat{\mathbf{u}}_{i+1}] = 	[\widetilde{\mathbf{u}}_{i+1,0}\mid \widetilde{\mathbf{u}}_{i+1,1}\mid\ldots\mid \widetilde{\mathbf{u}}_{i+1,i} \mid \widetilde{\mathbf{u}}_{i+1,i+1}].
		\end{align*}
		In other words,  $\widetilde{\mathbf{u}}_{i,j} = \widetilde{\mathbf{u}}_{i+1,j}$, $j=0,1\ldots, i$, and $ \widetilde{\mathbf{u}}_{i+1,i+1} = \mathbf{Q}_i\mathbf{Q}_i^\top \widehat{\mathbf{u}}_{i+1}$. Then, 
		\begin{align*}
			|\widetilde{\mathbf{u}}_{i+1,i+1} - \widehat{\mathbf{u}}_{i+1}|  = | \mathbf{Q}_i\mathbf{Q}_i^\top \widehat{\mathbf{u}}_{i+1} - \widehat{\mathbf{u}}_{i+1} | \le \textup{\texttt{tol}}. 
		\end{align*}
		Hence, $ \widetilde u_{h}^{i,j}=\widetilde{u}_{h}^{i+1,j} $ for $j=0,1,\ldots, i$, and 
		\begin{align*}
			\|\widehat u_{h}^{i+1}-\widetilde{u}_{h}^{i+1,i+1}  \|^2  &=  { (\widetilde{\mathbf{u}}_{i+1,i+1} - \widehat{\mathbf{u}}_{i+1})^\top \mathbf{M} (\widetilde{\mathbf{u}}_{i+1,i+1} - \widehat{\mathbf{u}}_{i+1})}\\
			&\le  {\sigma(\mathbf{M})} 	|\widetilde{\mathbf{u}}_{i+1,i+1} - \widehat{\mathbf{u}}_{i+1}|^2\\
			&\le \sigma(\mathbf{M})\textup{\texttt{tol}}^2\le Ch^2\textup{\texttt{tol}}^2.
		\end{align*}
		In  conclusion, we have 
		\begin{align*}
			\max\left\{\left\{\max_{0\le j \le i}		\|\widetilde  {u}_h^{i,j}-\widetilde{u}_h^{i+1, j}\|\right\},  \|\widehat {u}_h^{i+1}-\widetilde{u}_h^{i+1, i+1}\| \right\} \le  Ch\textup{\texttt{tol}}.
		\end{align*}
		
		{\bf Case 2:}  No truncation.  It is obvious that $\widetilde{\mathbf{u}}_{i,j} = \widetilde{\mathbf{u}}_{i+1,j}$ for  $j=0,1\ldots, i$, and $ \widetilde{\mathbf{u}}_{i+1,i+1} =   \widehat{\mathbf{u}}_{i+1}$.  Then  
		\begin{align*}
			\max\left\{\left\{\max_{0\le j \le i}		\|\widetilde  {u}_h^{i,j}-\widetilde{u}_h^{i+1, j}\|\right\},  \|\widehat {u}_h^{i+1}-\widetilde{u}_h^{i+1, i+1}\| \right\} \le Ch\textup{\texttt{tol}}.
		\end{align*}
		
		{\bf Case 3:} Singular value  truncation. By   \Cref{error_SingularValueTruncation} we have
		\begin{align*}
			|\widetilde{\mathbf{u}}_{i+1,j} - \widetilde{\mathbf{u}}_{i,j}|  \le \textup{\texttt{tol}}\; \textup{for}\; j=0,1,\ldots,i, \; \textup{and}\;	|\widetilde{\mathbf{u}}_{i+1,i+1} - \widehat{\mathbf{u}}_{i+1}|  \le \textup{\texttt{tol}}. 
		\end{align*}
		Following the same procedures as in case 1, we have
		\begin{align*}
			\|\widetilde  {u}_h^{i,j}-\widetilde{u}_h^{i+1, j}\|   \le Ch\textup{\texttt{tol}}\; \textup{for}\; j=0,1,\ldots,i, \; \textup{and}\;    \|\widehat {u}_h^{i+1}-\widetilde{u}_h^{i+1, i+1}\|  \le Ch \textup{\texttt{tol}}.
		\end{align*}
		Hence we have
		\begin{align*}
			\max\left\{\left\{\max_{0\le j \le i}		\|\widetilde  {u}_h^{i,j}-\widetilde{u}_h^{i+1, j}\|\right\},  \|\widehat {u}_h^{i+1}-\widetilde{u}_h^{i+1, i+1}\| \right\} \le Ch \textup{\texttt{tol}}.
		\end{align*}
	\end{proof}
\begin{lemma}\label{error_uhat_uh}
		At the $n$-th step, let $u_h^{n}$ and $\widehat{u}_h^{n}$ represent the solutions of \eqref{FEM} and \eqref{ISVD_eq1}, respectively. Throughout the entire process of the incremental SVD algorithm, the tolerance \textup{\texttt{tol}} is applied to both $p$-truncation and singular value truncation.  We establish the following error bound:
		\begin{align*}
			\left\|u_h^n-\widehat {u}_h^n \right\|  \le Chn \mathtt{tol},
		\end{align*} 
where $C$ is a positive constant independent of $h$. 
	\end{lemma}
	\begin{proof}
		We utilize mathematical induction to present our proof. Given that the initial condition remains consistent for both the finite element algorithm and the SVD algorithm, it follows that  $u_h^1 = \widehat u_h^1$. According to \Cref{error_uhat_uhtilde}, we can express the following inequality:
		\begin{align*}
			\max\left\{\|u_h^0 -  \widetilde u_h^{1, 0}\|, \|u_h^1 -  \widetilde u_h^{1, 1}\|\right\} \le Ch\textup{\texttt{tol}}.
		\end{align*}

		Assuming that for $i \le n-1$, the following holds:
		\begin{align}\label{Ass_i}
			\max_{0\le j\le i} \|u_h^j -  \widetilde u_h^{i, j}\|\le Chi \textup{\texttt{tol}}.
		\end{align}

		Next, we aim to prove the above inequality for $i+1$.  Using the triangle inequality, we get:
		\begin{align*}
			\max_{0\le j\le i+1} \|u_h^j -  &\widetilde u_h^{i+1, j}\|= \max\left\{ 	\max_{0\le j\le i} \|u_h^j -  \widetilde u_h^{i+1, j}\|, \|u_h^{i+1} -  \widetilde u_h^{i+1, i+1}\| \right\}\\
			&\le  \max\left\{ 	\max_{0\le j\le i}\left( \|u_h^j -  \widetilde u_h^{i, j}\|+ \| \widetilde u_h^{i, j}-\widetilde u_h^{i+1, j}\|\right), \|u_h^{i+1} -  \widetilde u_h^{i+1, i+1}\| \right\}.
		\end{align*} 
		By \eqref{Ass_i} and \Cref{error_uhat_uhtilde}, we have:
		\begin{align*}
			\max_{0\le j\le i}\left( \|u_h^j -  \widetilde u_h^{i, j}\|+ \| \widetilde u_h^{i, j}-\widetilde u_h^{i+1, j}\|\right) &\le \max_{0\le j\le i}  \|u_h^j -  \widetilde u_h^{i, j}\| +  \max_{0\le j\le i}   \| \widetilde u_h^{i, j}-\widetilde u_h^{i+1, j}\|\\
			& \le  Ch (i+1)\textup{\texttt{tol}}. 
		\end{align*}
		Also, by \eqref{Ass_i} and \Cref{error_uh_uhat}, we have:
		\begin{align*}
			\|u_h^{i+1} -  \widetilde u_h^{i+1, i+1}\| & \le 	\|u_h^{i+1} -  \widehat u_h^{i+1}\| + \|\widehat u_h^{i+1}-  \widetilde u_h^{i+1, i+1}\| \\
			&\le \max_{0\le j\le i} \| u_h^j-\widetilde{u}_h^{i,j} \| +  \|\widehat u_h^{i+1}-  \widetilde u_h^{i+1, i+1}\|\\
			& \le  Ch (i+1)\textup{\texttt{tol}}. 
		\end{align*}
		
		By selecting $i=n-2$, we obtain: 
		\begin{align*}
			\max_{0\le j\le n-1} \|u_h^j -  \widetilde u_h^{n-1, j}\| \le  Ch (n-1)\textup{\texttt{tol}}. 
		\end{align*}
		In accordance with \Cref{error_uh_uhat}  we deduce:
		\begin{align*}
			\left\|u_h^{n}-\widehat {u}_h^{n} \right\| \le 	\max_{0\le j\le n-1} \|u_h^j -  \widetilde u_h^{n-1, j}\| \le  Ch (n-1)\textup{\texttt{tol}}\le Chn\textup{\texttt{tol}}
		\end{align*}
where $C$ is a positive constant independent of $h$.
	\end{proof}
	
Hence, by applying the triangle inequality to \Cref{theorem401} and \Cref{error_uhat_uh}, we can derive the error estimate for $\|u(t_n)-\widehat{u}_{h}^{n} \| $. This completes the proof of \Cref{main_res}.
		
\section{Numerical experiments}
In this section, we use some examples to test the efficiency of  \Cref{fpdeandisvd} based on iFEM \cite{CheniFEM}. Since all three algorithms are for time discretization, we use the same spatial discretization for all numerical examples below and only record the computational time for the time marching process.
\subsection{Example 1}
{\color{blue} For the first example, we compare the numerical solutions produced by the direct method and the incremental SVD method at the final time $T=1$, reporting both the $L^2$-difference and the CPU times. }Let  $\Omega = (0,1)\times (0,1)$, we consider the equation \eqref{Fractional_PDEs} with
\begin{align*}
&\alpha = 0.5, ~~ ~~~T=1,  ~~~~u_0 = xy(1-x)(1-y),\\
 &f = 100\sin(2\pi t(x+y))x(1-x)y(1-y).	
\end{align*}
We use the linear finite element for the spatial discretization with different time steps $\tau$ for all three experiments below.  $h$ is the mesh size (max diameter of the triangles in the mesh). \textcolor{blue}{For the incremental SVD, we take \texttt{tol1}=\texttt{tol2}=\texttt{tol3} = $10^{-12}$ in \Cref{fpdeandisvd}.} For solving linear systems, we apply the Matlab built-in solver backslash ($\backslash$). We report the CPU time (seconds) \footnote{All the code for all examples in the paper has been created by the authors using Matlab R2025a and has been run on a MacBook Pro, M4 chip, 16 GB memory. We use the Matlab built-in function \texttt{tic}-\texttt{toc} to denote  the real simulation time.} of \Cref{fractionalPDEs} and \Cref{fpdeandisvd}  in \Cref{table_0,table_00}. We also compute the $L^2$-norm error between the solutions of \Cref{fractionalPDEs} and \Cref{fpdeandisvd}  at the final time $T=1$, we see that the error is close to machine precision. Furthermore,  \Cref{fpdeandisvd} is more efficient than \Cref{fractionalPDEs}.

\begin{table}[htbp]
	\centering
	\begingroup\color{blue}
	\caption{Comparison between the direct method and ISVD for \Cref{fractionalPDEs} with $\tau=10^{-3}$.}
	\label{table_0}
	\small
	\setlength{\tabcolsep}{8pt}
	\renewcommand{\arraystretch}{1.15}
	\begin{tabular}{cccc}
		\toprule
		$h$ & $\|u_h^n-\hat{u}_h^n\|_0$ & CPU time (s), ISVD & CPU time (s), Direct \\
		\midrule
		$1/2^5$ & $3\times10^{-14}$ & 0.24 & 0.37 \\
		$1/2^6$ & $3\times10^{-14}$ & 0.60 & 1.45 \\
		$1/2^7$ & $3\times10^{-14}$ & 2.54 & 7.38 \\
		$1/2^8$ & $3\times10^{-14}$ & 9.67 & 25.40 \\
		$1/2^9$ & $3\times10^{-14}$ & 53.53 & 126.72 \\
		\bottomrule
	\end{tabular}
	\endgroup
\end{table}

\begin{table}[htbp]
	\centering
	\begingroup\color{blue}
	\caption{Comparison between the direct method and ISVD for \Cref{fractionalPDEs} with $\tau=10^{-4}$.}
	\label{table_00}
	\small
	\setlength{\tabcolsep}{8pt}
	\renewcommand{\arraystretch}{1.15}
	\begin{tabular}{cccc}
		\toprule
		$h$ & $\|u_h^n-\hat{u}_h^n\|_0$ & CPU time (s), ISVD & CPU time (s), Direct \\
		\midrule
		$1/2^5$ & $3\times10^{-14}$ & 1.56 & 2.90 \\
		$1/2^6$ & $3\times10^{-14}$ & 2.26 & 6.10 \\
		$1/2^7$ & $3\times10^{-14}$ & 4.56 & 23.45 \\
		$1/2^8$ & $3\times10^{-14}$ & 15.70 & 93.63 \\
		$1/2^9$ & $3\times10^{-14}$ & 56.97 & 368.13 \\
		\bottomrule
	\end{tabular}
	\endgroup
\end{table}

\subsection{Example 2}
\begin{figure}[H]
    \centering
    \includegraphics[width=0.95\textwidth]{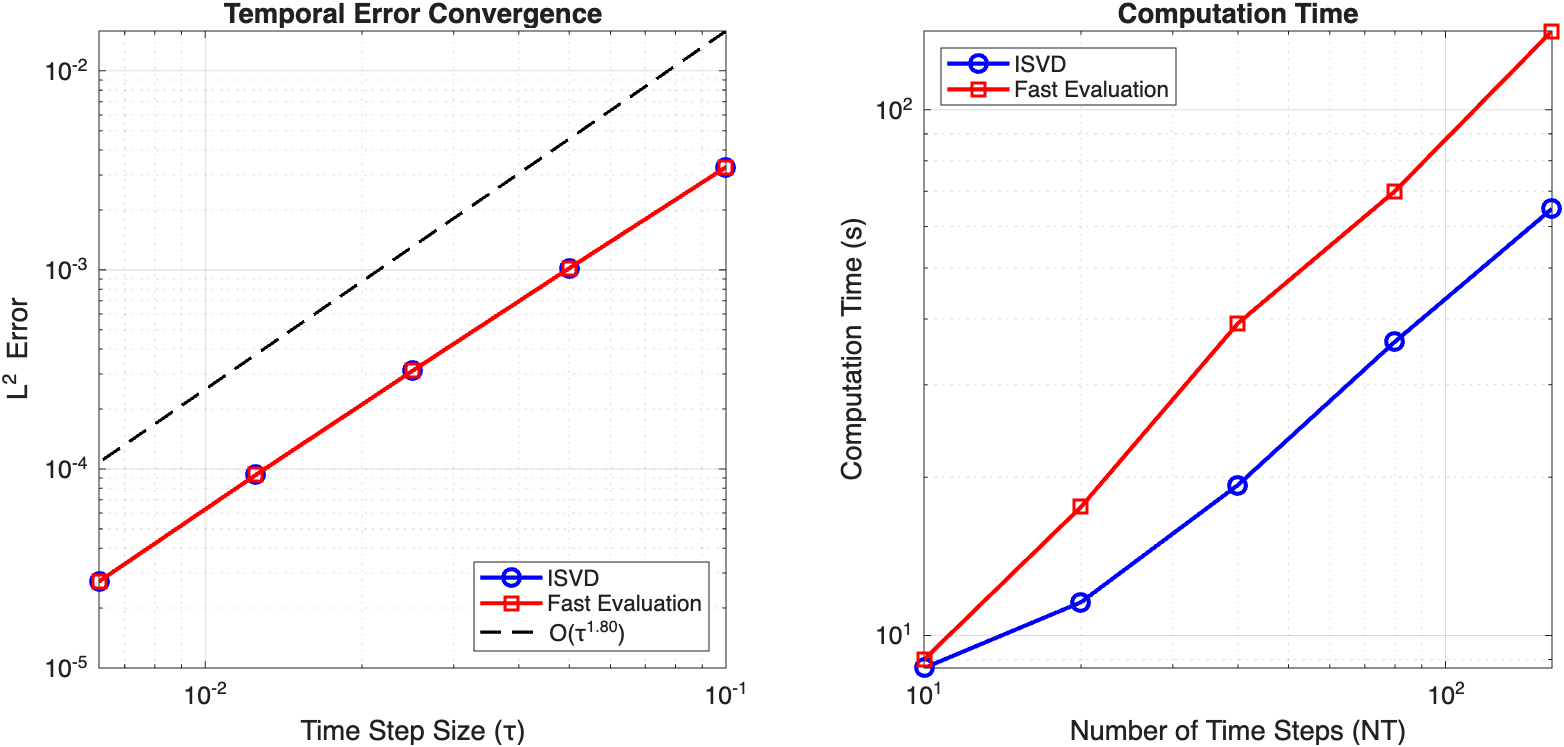}
    \caption{Temporal convergence rate and comparison of ISVD and fast evaluation methods for sub-diffusion equation with $\alpha = 0.2,T=1$}
    \label{fig:alp02_time}
\end{figure}

\begin{figure}[H]
    \centering
        \includegraphics[width=0.95\textwidth]{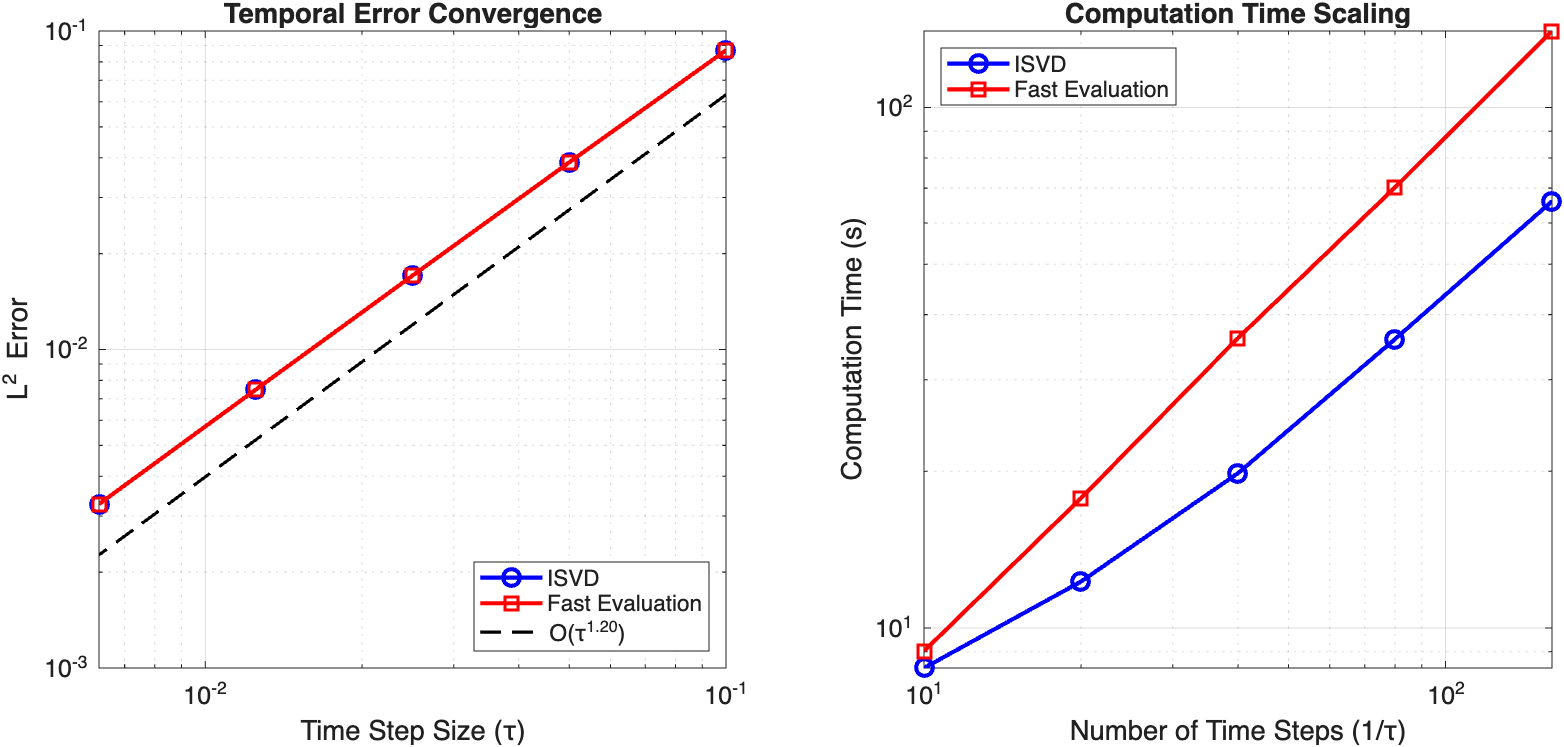}
    \caption{Temporal convergence rate and comparison of ISVD and fast evaluation methods for sub-diffusion equation with $\alpha = 0.8,T=1$}
    \label{fig:alp08_time}
\end{figure}

\begin{figure}[H]
    \centering
        \includegraphics[width=0.95\textwidth]{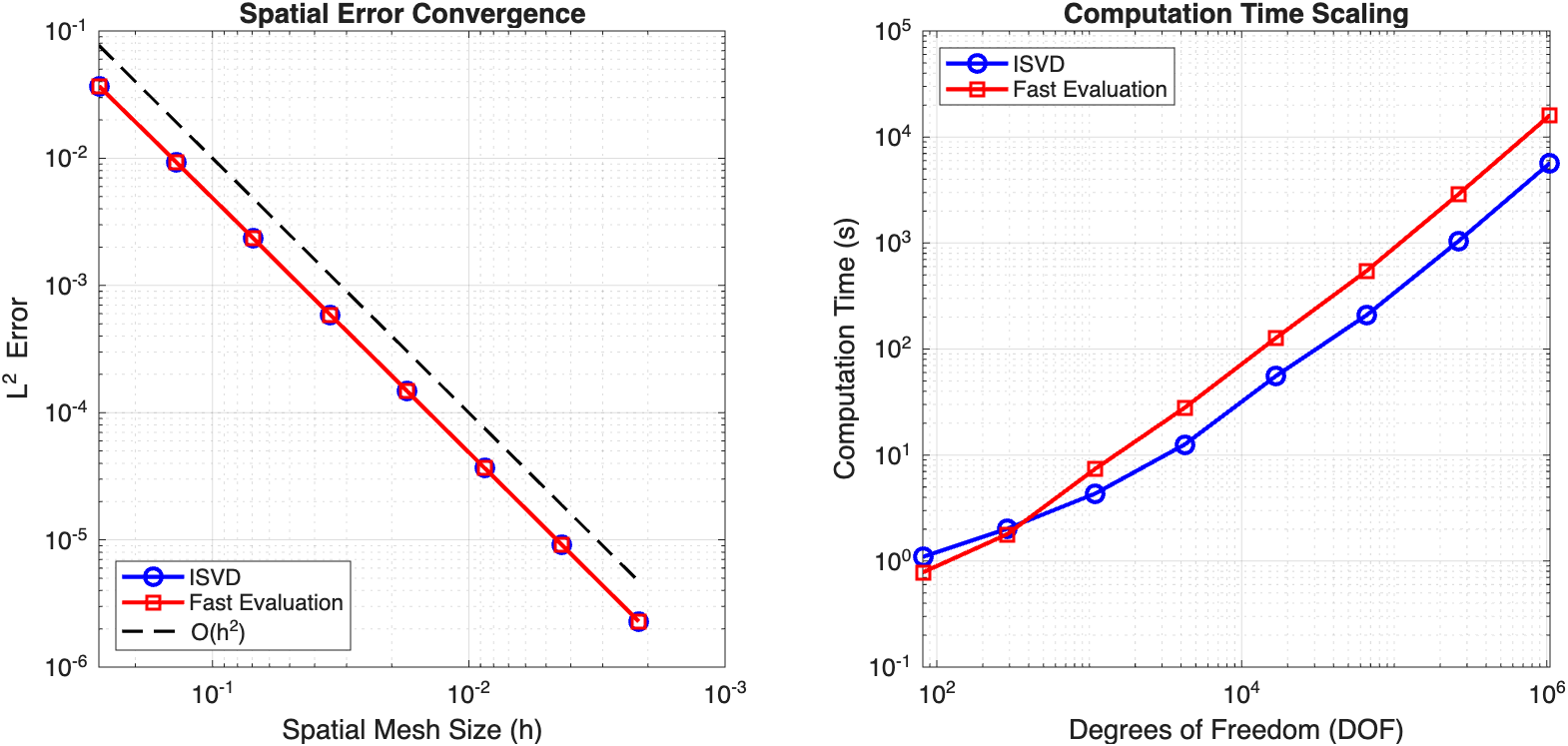}
    \caption{Spatial convergence rate and comparison of ISVD and fast evaluation methods for sub-diffusion equation with $\alpha = 0.2,T=1$ on refined mesh.}
    \label{fig:alp02_space}
\end{figure}

{\color{blue} To evaluate the computational efficiency and accuracy of the incremental SVD scheme, we use the example introduced in \cite{Zhuang2019variableCoeff, Ren2018fast}. We choose the fast evaluation method in \cite{Zhang_fast_evaluation2017} as the benchmark because it is a representative fast-history algorithm within the same $L1$/finite element framework, which allows a fair comparison under identical spatial and temporal discretizations.}  This canonical problem examines a variable-coefficient fractional sub-diffusion equation:
\begin{equation}\notag
\left\{
\begin{aligned}
\partial_{t}^{\alpha} u(\mathbf{x},t) - \nabla\cdot\left(b(\mathbf{x})\nabla u(\mathbf{x},t)\right) &= f(\mathbf{x},t) && (\mathbf{x},t) \in \Omega \times (0,T], \\
u(\mathbf{x},t) &= 0 && (\mathbf{x},t) \in \partial \Omega \times (0,T], \\
u(\mathbf{x},0) &= u_0(\mathbf{x}) && \mathbf{x} \in \Omega.
\end{aligned}
\right.
\end{equation}
Set $\Omega=(0,\pi)^2$, $b(\mathbf{x})=\sin x\sin y+0.1$ with right-hand side
 \[
 f=\frac{\Gamma(3+\alpha)}{2}t^2\sin x\sin y-t^{2+\alpha}\left(\sin^2x+\sin^2y-4\sin^2x\sin^2y-0.2\sin x\sin y\right).
 \]
 The analytic solution is given by $u = t^{2+\alpha}\sin x\sin y.$  Figures \ref{fig:alp02_time},\ref{fig:alp08_time},\ref{fig:alp02_space} demonstrate that both ISVD and Fast Evaluation effectively simulate the physical problem. When fixing the mesh size at $h=\frac{\pi}{1000}$, both methods achieve a $O(\tau^{2-\alpha})$ convergence rate, as shown in Figures \ref{fig:alp02_time} and \ref{fig:alp08_time}. {\color{blue}For this fixed mesh size, the CPU-time comparison shows that ISVD remains competitive with the fast evaluation method across the tested time-step sizes while producing essentially identical errors. Figure~\ref{fig:alp02_space} shows the spatial convergence behavior and CPU times for both algorithms when $\tau=10^{-4}$, the fast-method tolerance is $\varepsilon=10^{-9}$, and the ISVD tolerances are $\mathtt{tol1}=\mathtt{tol2}=\mathtt{tol3}=10^{-12}$. The detailed results in Tables~\ref{tab: combined_results_vertical1}, \ref{tab: combined_results_vertical2}, and \ref{tab: combined_results_vertical3} show that both approaches preserve the same accuracy, while ISVD becomes faster on the refined meshes considered here. The omitted direct-method entries in Tables~\ref{tab: combined_results_vertical2} and \ref{tab: combined_results_vertical3} correspond to parameter settings for which the required history storage exceeded the available 16~GB memory on our test machine.}

\begin{sidewaystable}
\centering
\color{blue}
\arrayrulecolor{blue} 
\caption{Spatial convergence rates and computational time with $\tau = 10^{-4},T=1, \varepsilon=10^{-9},\mathtt{tol} = 10^{-12}$}
\label{tab: combined_results_vertical1}
\small
\setlength{\tabcolsep}{4pt}
\begin{tabular}{@{}c
  S[table-format=1.2e-1] S[table-format=1.2] S[table-format=4.2]
  S[table-format=1.2e-1] S[table-format=1.2] S[table-format=4.2]
  S[table-format=1.2e-1] S[table-format=1.2] S[table-format=4.2]@{}}
\toprule
\multicolumn{10}{c}{$\alpha = 0.2$} \\
\cmidrule(lr){1-10}
\multirow{2}{*}{$h$} & \multicolumn{3}{c}{ISVD} & \multicolumn{3}{c}{Fast Evaluation} & \multicolumn{3}{c}{Direct Method} \\
\cmidrule(lr){2-4} \cmidrule(lr){5-7} \cmidrule(lr){8-10}
& {Error} & {Rate} & {CPU(s)\footnote{All errors are in $L^2$ norm. Computational times in seconds. Convergence rates (Rate) are calculated between consecutive mesh refinements.}}
& {Error} & {Rate} & {CPU(s)} & {Error} & {Rate} & {CPU(s)} \\
\midrule
$\pi/20$  & 5.99e-3 &      &  2.71  & 5.99e-3 &      &   2.78 & 5.99e-3 &      &  10.39 \\
$\pi/40$  & 1.50e-3 & 1.99 &  5.49  & 1.50e-3 & 1.99 &  10.73 & 1.50e-3 & 1.99 &  41.68 \\
$\pi/80$  & 3.76e-4 & 2.00 & 18.39  & 3.76e-4 & 2.00 &  43.25 & 3.76e-4 & 2.00 & 120.54 \\
$\pi/160$ & 9.40e-5 & 2.00 & 80.63  & 9.40e-5 & 2.00 & 196.15 & 9.40e-5 & 2.00 & 450.99 \\
$\pi/320$ & 2.35e-5 & 2.00 & 363.10 & 2.35e-5 & 2.00 & 943.95 & 2.35e-5 & 2.00 & 1552.71 \\
\addlinespace[2mm]
\multicolumn{10}{c}{$\alpha = 0.5$} \\
\cmidrule(lr){1-10}
\multirow{2}{*}{$h$} & \multicolumn{3}{c}{ISVD} & \multicolumn{3}{c}{Fast Evaluation} & \multicolumn{3}{c}{Direct Method} \\
\cmidrule(lr){2-4} \cmidrule(lr){5-7} \cmidrule(lr){8-10}
& {Error} & {Rate} & {CPU(s)} & {Error} & {Rate} & {CPU(s)} & {Error} & {Rate} & {CPU(s)} \\
\midrule
$\pi/20$  & 5.32e-3 &      &  2.91  & 5.32e-3 &      &   2.98 & 5.32e-3 &      &   9.40 \\
$\pi/40$  & 1.33e-3 & 1.99 &  5.67  & 1.33e-3 & 1.99 &  11.32 & 1.33e-3 & 1.99 &  38.66 \\
$\pi/80$  & 3.33e-4 & 2.00 & 18.26  & 3.33e-4 & 2.00 &  42.80 & 3.33e-4 & 2.00 & 116.08 \\
$\pi/160$ & 8.30e-5 & 2.00 & 80.58  & 8.30e-5 & 2.00 & 193.46 & 8.30e-5 & 2.00 & 450.49 \\
$\pi/320$ & 2.05e-5 & 2.00 & 372.73 & 2.05e-5 & 2.00 & 928.61 & 2.05e-5 & 2.00 & 1540.99 \\
\addlinespace[2mm]
\multicolumn{10}{c}{$\alpha = 0.8$} \\
\cmidrule(lr){1-10}
\multirow{2}{*}{$h$} & \multicolumn{3}{c}{ISVD} & \multicolumn{3}{c}{Fast Evaluation} & \multicolumn{3}{c}{Direct Method} \\
\cmidrule(lr){2-4} \cmidrule(lr){5-7} \cmidrule(lr){8-10}
& {Error} & {Rate} & {CPU(s)} & {Error} & {Rate} & {CPU(s)} & {Error} & {Rate} & {CPU(s)} \\
\midrule
$\pi/20$  & 4.56e-3 &      &  4.50  & 4.56e-3 &      &   2.84 & 4.56e-3 &      &   9.44 \\
$\pi/40$  & 1.13e-3 & 2.01 &  9.39  & 1.13e-3 & 2.01 &  10.98 & 1.13e-3 & 2.01 &  41.32 \\
$\pi/80$  & 2.77e-4 & 2.03 & 30.53  & 2.77e-4 & 2.03 &  46.36 & 2.77e-4 & 2.03 & 123.43 \\
$\pi/160$ & 6.65e-5 & 2.06 & 101.44 & 6.65e-5 & 2.06 & 196.52 & 6.65e-5 & 2.06 & 446.31 \\
$\pi/320$ & 1.66e-5 & 2.00 & 434.08 & 1.66e-5 & 2.00 & 953.74 & 1.66e-5 & 2.00 & 1393.94 \\
\bottomrule
\end{tabular}
\arrayrulecolor{black} 
\end{sidewaystable}
\begin{sidewaystable}
\centering
\color{blue}
\arrayrulecolor{blue} 
\caption{Spatial convergence rates and computational performance with $\tau = 10^{-4},T=5, \varepsilon=10^{-9},\mathtt{tol} = 10^{-12}$}
\label{tab: combined_results_vertical2}
\small
\setlength{\tabcolsep}{4pt}
\begin{tabular}{@{}c
  S[table-format=1.2e-1] S[table-format=1.2] S[table-format=4.2]
  S[table-format=1.2e-1] S[table-format=1.2] S[table-format=4.2]
  S[table-format=1.2e-1] S[table-format=1.2] S[table-format=4.2]@{}}
\toprule
\multicolumn{10}{c}{$\alpha = 0.2$} \\
\cmidrule(lr){1-10}
\multirow{2}{*}{$h$} & \multicolumn{3}{c}{ISVD} & \multicolumn{3}{c}{Fast Evaluation} & \multicolumn{3}{c}{Direct Method} \\
\cmidrule(lr){2-4} \cmidrule(lr){5-7} \cmidrule(lr){8-10}
& {Error} & {Rate} & {CPU(s)\footnote{All errors are in $L^2$ norm. Computational times in seconds. Convergence rates (Rate) are calculated between consecutive mesh refinements.}}
& {Error} & {Rate} & {CPU(s)} & {Error} & {Rate} & {CPU(s)} \\
\midrule
$\pi/20$  &2.28e-01 &      &19.38    &2.28e-01  &      &13.96    &2.28e-01  &      &194.16   \\
$\pi/40$  &5.72e-02 & 1.99 &33.58     &5.72e-02  &1.99  &55.00   &5.72e-02  &1.99  &938.36   \\
$\pi/80$  &1.43e-02 & 1.99 &95.26   &1.43e-02  &1.99  &212.31   &1.43e-02  &1.99  &2597.31  \\
$\pi/160$ &3.58e-03  &2.00  &395.84    &3.58e-03  &2.00  &955.88  &3.58e-03  &2.00  &11141.63   \\
$\pi/320$ &  8.95e-04 &2.00  &1861.44  &8.95e-04  & 2.00 &4458.69  &\multicolumn{1}{c}{--}  &\multicolumn{1}{c}{--}  &\multicolumn{1}{c}{--}  \\
\addlinespace[2mm]
\multicolumn{10}{c}{$\alpha = 0.5$} \\
\cmidrule(lr){1-10}
\multirow{2}{*}{$h$} & \multicolumn{3}{c}{ISVD} & \multicolumn{3}{c}{Fast Evaluation} & \multicolumn{3}{c}{Direct Method} \\
\cmidrule(lr){2-4} \cmidrule(lr){5-7} \cmidrule(lr){8-10}
& {Error} & {Rate} & {CPU(s)} & {Error} & {Rate} & {CPU(s)} & {Error} & {Rate} & {CPU(s)} \\
\midrule
$\pi/20$  & 3.85e-1 &      &29.51    & 3.85e-1 &      &14.14    &3.85e-01  &      &196.65    \\
$\pi/40$  & 9.65e-2 & 1.99 &  50.21  & 9.65e-2 & 1.99 &  58.06 &9.65e-02  &1.99  &937.11    \\
$\pi/80$  & 2.41e-2 & 2.00 &  134.38 & 2.41e-2 & 2.00 &  208.96 &2.41e-02  &2.00  &2588.03  \\
$\pi/160$ & 6.03e-3 & 2.00 &  443.73 & 6.03e-3 & 2.00 & 1102.06 &6.03e-03  &2.00  &10822.21  \\
$\pi/320$ & 1.51e-03 &2.00  &1999.21  &1.51e-03  & 2.00 &3672.00  &\multicolumn{1}{c}{--}  &\multicolumn{1}{c}{--}  &\multicolumn{1}{c}{--}  \\
\addlinespace[2mm]
\multicolumn{10}{c}{$\alpha = 0.8$} \\
\cmidrule(lr){1-10}
\multirow{2}{*}{$h$} & \multicolumn{3}{c}{ISVD} & \multicolumn{3}{c}{Fast Evaluation} & \multicolumn{3}{c}{Direct Method} \\
\cmidrule(lr){2-4} \cmidrule(lr){5-7} \cmidrule(lr){8-10}
& {Error} & {Rate} & {CPU(s)} & {Error} & {Rate} & {CPU(s)} & {Error} & {Rate} & {CPU(s)} \\
\midrule
$\pi/20$  &6.38e-1  &      &31.21   & 6.384e-1 &      &14.27    &6.38e-01  &      &196.42    \\
$\pi/40$  & 1.60e-01 &1.99  &  59.02    & 1.60e-01 & 1.99 &56.35    & 1.60e-01 &1.99  &945.57   \\
$\pi/80$  & 4.00e-02  &1.99  &   157.21&4.00e-02  & 2.00 &201.68   &4.00e-02  &2.00  &2589.34  \\
$\pi/160$ &9.92e-03   & 2.00 & 470.29 &9.92e-03  &2.00&937.05  &9.92e-03  &2.00  &\multicolumn{1}{c}{10653.82}  \\
$\pi/320$ & 2.41e-03 &2.00  & 2091.78  &2.41e-03 & 2.00 &4390.92  &  \multicolumn{1}{c}{--} &\multicolumn{1}{c}{--}  &\multicolumn{1}{c}{--}  \\
\bottomrule
\end{tabular}
\arrayrulecolor{black} 
\end{sidewaystable}
\begin{sidewaystable}
\centering
\color{blue}
\arrayrulecolor{blue} 
\caption{Spatial convergence rates and computational performance with $\tau = 5\times10^{-5},T=10, \varepsilon=10^{-9},\mathtt{tol1} = 10^{-10},\mathtt{tol2} = 10^{-14},\mathtt{tol3} = 10^{-8}$}
\label{tab: combined_results_vertical3}
\small
\setlength{\tabcolsep}{4pt}
\begin{tabular}{@{}c
  S[table-format=1.2e-1] S[table-format=1.2] S[table-format=4.2]
  S[table-format=1.2e-1] S[table-format=1.2] S[table-format=4.2]
  S[table-format=1.2e-1] S[table-format=1.2] S[table-format=4.2]@{}}
\toprule
\multicolumn{10}{c}{$\alpha = 0.2$} \\
\cmidrule(lr){1-10}
\multirow{2}{*}{$h$} & \multicolumn{3}{c}{ISVD} & \multicolumn{3}{c}{Fast Evaluation} & \multicolumn{3}{c}{Direct Method} \\
\cmidrule(lr){2-4} \cmidrule(lr){5-7} \cmidrule(lr){8-10}
& {Error} & {Rate} & {CPU(s)\footnote{All errors are in $L^2$ norm. Computational times in seconds. Convergence rates (Rate) are calculated between consecutive mesh refinements.}}
& {Error} & {Rate} & {CPU(s)} & {Error} & {Rate} & {CPU(s)} \\
\midrule
$\pi/20$  &\multicolumn{1}{c}{1.09} &  \multicolumn{1}{c}{--}    &220.43    &\multicolumn{1}{c}{1.09}   & \multicolumn{1}{c}{--}    &27.15    &\multicolumn{1}{c}{1.09}  &      &754.92   \\
$\pi/40$  &2.74e-01  &1.99   &329.98     &2.74e-01  &1.99  &106.15   &2.74e-01  &1.99  &3651.54   \\
$\pi/80$  &6.84e-02  &2.00  & 653.28  &6.84e-02  &2.00  &415.04   &6.84e-02  & 2.00 &10187.13  \\
$\pi/160$ &1.71e-02  &2.00  &1855.10  &1.71e-02  &2.00  &1716.88  &  \multicolumn{1}{c}{--} & \multicolumn{1}{c}{--}   & \multicolumn{1}{c}{--}    \\
$\pi/320$ & 4.28e-03  &2.00 & 7463.04 &4.28e-03  &2.00  &7476.68  & \multicolumn{1}{c}{--}   & \multicolumn{1}{c}{--}   &   \multicolumn{1}{c}{--} \\
\addlinespace[2mm]
\multicolumn{10}{c}{$\alpha = 0.5$} \\
\cmidrule(lr){1-10}
\multirow{2}{*}{$h$} & \multicolumn{3}{c}{ISVD} & \multicolumn{3}{c}{Fast Evaluation} & \multicolumn{3}{c}{Direct Method} \\
\cmidrule(lr){2-4} \cmidrule(lr){5-7} \cmidrule(lr){8-10}
& {Error} & {Rate} & {CPU(s)} & {Error} & {Rate} & {CPU(s)} & {Error} & {Rate} & {CPU(s)} \\
\midrule
$\pi/20$  &\multicolumn{1}{c}{2.39} &\multicolumn{1}{c}{--}       &487.01   & \multicolumn{1}{c}{2.38}  &      &28.65    &\multicolumn{1}{c}{2.38}   &      &765.57    \\
$\pi/40$  &5.98e-01 & 1.99  &   471.91 &5.98e-01  &  1.99&110.39   &  5.98e-01& 1.99 & 3707.56   \\
$\pi/80$  &1.50e-01  &2.00  &794.54   &  1.50e-01& 2.00 &428.58   &1.50e-01  & 2.00 &10717.13  \\
$\pi/160$ &3.74e-02  &2.00  &  2013.47 &3.74e-02  &2.00  &1712.27  & \multicolumn{1}{c}{--}   & \multicolumn{1}{c}{--}   & \multicolumn{1}{c}{--}   \\
$\pi/320$ &9.35e-03  & 2.00 & 7714.61 &9.35e-03  &2.00  &7668.15  &  \multicolumn{1}{c}{--}  & \multicolumn{1}{c}{--}   &  \multicolumn{1}{c}{--}  \\
\addlinespace[2mm]
\multicolumn{10}{c}{$\alpha = 0.8$} \\
\cmidrule(lr){1-10}
\multirow{2}{*}{$h$} & \multicolumn{3}{c}{ISVD} & \multicolumn{3}{c}{Fast Evaluation} & \multicolumn{3}{c}{Direct Method} \\
\cmidrule(lr){2-4} \cmidrule(lr){5-7} \cmidrule(lr){8-10}
& {Error} & {Rate} & {CPU(s)} & {Error} & {Rate} & {CPU(s)} & {Error} & {Rate} & {CPU(s)} \\
\midrule
$\pi/20$  &\multicolumn{1}{c}{5.11}   & \multicolumn{1}{c}{--}      &1864.22     &\multicolumn{1}{c}{5.11}   &      & 28.22   &\multicolumn{1}{c}{5.11}   &      &756.05    \\
$\pi/40$  &  \multicolumn{1}{c}{1.28} &2.00  & 1908.85     &\multicolumn{1}{c}{1.28}   & 2.00 &109.00    & \multicolumn{1}{c}{1.28}  &  2.00&3669.58   \\
$\pi/80$  & 3.20e-01  &2.00   &1970.05   &3.20e-01 &1.99  &428.48   &3.20e-01  &2.00  &\multicolumn{1}{c}{11043.23}  \\
$\pi/160$ & 7.98e-02  &2.00  &  3099.71&7.98e-02  &2.00&1891.93  & \multicolumn{1}{c}{--}   &   \multicolumn{1}{c}{--} & \multicolumn{1}{c}{--}   \\
$\pi/320$ & 1.98e-02 &2.00  & 7019.23 &1.98e-02 &2.01  &7233.63  &  \multicolumn{1}{c}{--}   &   \multicolumn{1}{c}{--} &  \multicolumn{1}{c}{--}  \\
\bottomrule
\end{tabular}
\arrayrulecolor{black} 
\end{sidewaystable}
{\color{blue}Tables~\ref{tab: combined_results_vertical1}, \ref{tab: combined_results_vertical2}, and \ref{tab: combined_results_vertical3} make two additional points clear. First, for $\tau = 10^{-4}$ and $T=1$, ISVD and the fast evaluation method achieve essentially the same level of accuracy, while ISVD becomes faster on the finer meshes. Second, the longer-time runs with $T=5$ and $T=10$ show why history compression is needed: as the number of time levels grows, the direct method quickly becomes impractical on the finer meshes, whereas the two compressed-history approaches remain feasible.}
{\color{blue}This behavior indicates that ISVD is effective for simulations with fine spatial discretizations, and a similar trend is observed in the three-dimensional experiment reported in Example~3.}
\begin{table}[H]
\centering
\color{blue}
\arrayrulecolor{blue} 
\caption{Impact of approximation tolerance ($\varepsilon$) on the fast evaluation method \cite{Zhang_fast_evaluation2017} ($\tau  = 10^{-3}$, $N_{dof}=10201$).}
\label{tab:fast_epsilon_impact1}
\renewcommand{\arraystretch}{1.2}
\begin{tabular}{@{}c cc c cc@{}}
\toprule
\textbf{Case} & \textbf{ $\varepsilon$} & \textbf{$N_{\text{exp}}$} & \textbf{$L^2$ Error} & \textbf{Total Time (s)} & \textbf{Time/Step (s)} \\ 
\midrule
1 & $10^{-3}$  & 20 & $2.1947 \times 10^{-4}$ & 6.68 & $6.6832 \times 10^{-3}$ \\
2 & $10^{-6}$  & 35 & $2.0404 \times 10^{-4}$ & 6.00 & $6.0047 \times 10^{-3}$ \\
3 & $10^{-9}$  & 43 & $2.0403 \times 10^{-4}$ & 6.37 & $6.3719 \times 10^{-3}$ \\
4 & $10^{-10}$ & 53 & $2.0403 \times 10^{-4}$ & 6.29 & $6.2919 \times 10^{-3}$ \\
5 & $10^{-12}$ & 62 & $2.0403 \times 10^{-4}$ & 6.21 & $6.2081 \times 10^{-3}$ \\
\bottomrule
\end{tabular}
\arrayrulecolor{black} 
\end{table}
\begin{table}[H]
\centering
\color{blue}
\arrayrulecolor{blue} 
\caption{Impact of approximation tolerance ($\varepsilon$) on the fast evaluation method($\tau  = 10^{-3}$, $N_{dof}=1002001$).}
\label{tab:fast_epsilon_impact2}
\renewcommand{\arraystretch}{1.2}
\begin{tabular}{@{}l cc c cc@{}}
\toprule
\textbf{Test} & \textbf{$\varepsilon$} & \textbf{$N_{\text{exp}}$} & \textbf{$L^2$ Error} & \textbf{Total Time (s)} & \textbf{Time/Step (s)} \\ 
\midrule
1 & $10^{-3}$  & 20 & $1.4160 \times 10^{-5}$ & 964.62  & $9.6462 \times 10^{-1}$ \\
2 & $10^{-6}$  & 35 & $2.4137 \times 10^{-5}$ & 884.58  & $8.8458 \times 10^{-1}$ \\
3 & $10^{-9}$  & 43 & $2.4165 \times 10^{-5}$ & 1019.00 & $1.0190 \times 10^{0}$  \\
4 & $10^{-10}$ & 53 & $2.4165 \times 10^{-5}$ & 1252.42 & $1.2524 \times 10^{0}$  \\
5 & $10^{-12}$ & 62 & $2.4165 \times 10^{-5}$ & 911.20  & $9.1120 \times 10^{-1}$ \\
\bottomrule
\end{tabular}
\arrayrulecolor{black} 
\end{table}
\textcolor{blue}{To make the fast-method tuning explicit, Tables~\ref{tab:fast_epsilon_impact1} and~\ref{tab:fast_epsilon_impact2} report the approximation tolerance $\varepsilon$, the resulting number of exponentials $N_{\exp}$, the final $L^2$ error, and the CPU time. In both tests, the error decreases when $\varepsilon$ is tightened from $10^{-3}$ to $10^{-9}$ and then essentially saturates, while $N_{\exp}$ increases from $20$ to $43$ and then to $62$ for stricter tolerances. We therefore use $\varepsilon=10^{-9}$ in the main comparisons, because it reaches the discretization-limited accuracy without introducing unnecessary exponential terms.}

\begin{table}[H]
\centering
\color{blue}
\arrayrulecolor{blue} 
\caption{Impact of truncation and orthogonality tolerances $\mathtt{tol1,tol2,tol3}$ on the ISVD method ($N_t = T/ \tau = 1000$, $N_{dof}=10201, \alpha=0.8$).}
\label{tab:isvd_tolerances1}
\renewcommand{\arraystretch}{1.2} 
\begin{tabular}{@{}c ccc c ccc@{}}
\toprule
\textbf{Case} & \textbf{tol1} & \textbf{tol2} & \textbf{tol3} & \textbf{Rank} & \textbf{ Error} & \textbf{Total Time (s)} & \textbf{Time/Step (s)} \\ 
\midrule
1  & $10^{-12}$ & $10^{-14}$ & $10^{-10}$ & 7 & $3.5064\times 10^{-4}$ & 10.08 & $1.0076\times 10^{-2}$ \\
2 & $10^{-6}$  & $10^{-14}$ & $10^{-10}$ & 6 & $1.1339 \times 10^{-3}$ & 10.14 & $1.0142\times10^{-2}$ \\
3 & $10^{-12}$ & $10^{-13}$  & $10^{-10}$ & 7 & $3.5064 \times 10^{-4}$ & 9.48 & $9.4819\times10^{-3}$ \\
4 & $10^{-12}$ & $10^{-14}$ & $10^{-4}$  & 1 & $7.4808 \times 10^{-4}$ & 10.26 & $1.0257\times10^{-2}$ \\ 
5 & $10^{-12}$ & $10^{-12} $ &$10^{-12}$ &16  &$ 3.5065\times 10^{-4}$   & 10.41   & $ 1.0412\times10^{-2}$\\
\bottomrule
\end{tabular}
\arrayrulecolor{black} 
\end{table}

\begin{table}[H]
\centering
\color{blue}
\arrayrulecolor{blue} 
\caption{Impact of truncation and orthogonality tolerances $\mathtt{tol1,tol2,tol3}$ on the ISVD method ($N_t = T/\tau = 1000$, $N_{dof}=40401$, $\alpha=0.8$).}
\label{tab:isvd_tolerances2}
\renewcommand{\arraystretch}{1.2} 
\begin{tabular}{@{}c ccc c ccc@{}}
\toprule
\textbf{Case} & \textbf{tol1} & \textbf{tol2} & \textbf{tol3} & \textbf{Rank} & \textbf{Error} & \textbf{Total Time (s)} & \textbf{Time/Step (s)} \\ 
\midrule
1 & $10^{-12}$ & $10^{-14}$ & $10^{-10}$ & 5  & $3.6054 \times 10^{-4}$ & 45.95 & $4.5947 \times 10^{-2}$ \\
2 & $10^{-6}$  & $10^{-14}$ & $10^{-10}$ & 4  & $2.5613 \times 10^{-3}$ & 45.99 & $4.5988 \times 10^{-2}$ \\
3 & $10^{-12}$ & $10^{-13}$ & $10^{-10}$ & 5  & $3.6054 \times 10^{-4}$ & 45.19 & $4.5185 \times 10^{-2}$ \\
4 & $10^{-12}$ & $10^{-14}$ & $10^{-4}$  & 5  & $7.5871 \times 10^{-4}$ & 46.25 & $4.6248 \times 10^{-2}$ \\
5 & $10^{-12}$ & $10^{-12}$ & $10^{-12}$ & 16 & $3.6054 \times 10^{-4}$ & 48.88 & $4.8875 \times 10^{-2}$ \\
\bottomrule
\end{tabular}
\arrayrulecolor{black} 
\end{table}
\begin{table}[H]
\centering
\color{blue}
\arrayrulecolor{blue} 
\caption{Impact of truncation and orthogonality tolerances $\mathtt{tol1,tol2,tol3}$ on the ISVD method ($N_t = T/\tau = 1000$, $N_{dof}=160801$, $\alpha=0.8$).}
\label{tab:isvd_tolerances3}
\renewcommand{\arraystretch}{1.2} 
\begin{tabular}{@{}c ccc c ccc@{}}
\toprule
\textbf{Case} & \textbf{tol1} & \textbf{tol2} & \textbf{tol3} & \textbf{Rank} & \textbf{ Error} & \textbf{Total Time (s)} & \textbf{Time/Step (s)} \\ 
\midrule
1 & $10^{-12}$ & $10^{-14}$ & $10^{-10}$ & 5  & $3.6351 \times 10^{-4}$ & 239.71 & $2.3971 \times 10^{-1}$ \\
2 & $10^{-6}$  & $10^{-14}$ & $10^{-10}$ & 1  & $3.7952 \times 10^{-3}$ & 220.86 & $2.2086 \times 10^{-1}$ \\
3 & $10^{-12}$ & $10^{-13}$ & $10^{-10}$ & 5  & $3.6351 \times 10^{-4}$ & 231.25 & $2.3125 \times 10^{-1}$ \\
4 & $10^{-12}$ & $10^{-14}$ & $10^{-4}$  & 1  & $7.6778 \times 10^{-4}$ & 215.88 & $2.1588 \times 10^{-1}$ \\
5 & $10^{-12}$ & $10^{-12}$ & $10^{-12}$ & 15 & $3.6347 \times 10^{-4}$ & 232.72 & $2.3272 \times 10^{-1}$ \\
\bottomrule
\end{tabular}
\arrayrulecolor{black} 
\end{table}
{\color{blue}Tables~\ref{tab:isvd_tolerances1}, \ref{tab:isvd_tolerances2}, and \ref{tab:isvd_tolerances3} also clarify the role of the ISVD tolerances. Relaxing $\mathtt{tol1}$ or $\mathtt{tol3}$ reduces the retained rank and can slightly reduce the CPU time, but it also increases the final error; by contrast, changing $\mathtt{tol2}$ within the tested near-machine-precision range has almost no visible effect, which is consistent with its role as a reorthogonalization safeguard. In all cases the retained rank remains very small compared with both $N_{dof}$ and $N_t$, confirming the low-rank assumption used by the algorithm.

This also makes the storage comparison concrete. For the representative case $N_{dof}=160801$ and $N_t=1000$, Table~\ref{tab:isvd_tolerances3} gives retained rank $r=15$ for the strict choice $\mathtt{tol1}=\mathtt{tol2}=\mathtt{tol3}=10^{-12}$. The dominant ISVD storage is then approximately $(N_{dof}+N_t)r+r^2 = 2.43\times 10^6$ doubles, whereas the direct method requires $N_{dof}N_t = 1.61\times 10^8$ doubles. For the fast method, Tables~\ref{tab:fast_epsilon_impact1}--\ref{tab:fast_epsilon_impact2} show that $\varepsilon=10^{-9}$ gives $N_{\exp}=43$, corresponding to about $N_{dof}N_{\exp} = 6.91\times 10^6$ doubles for the same $N_{dof}$. Thus, in this tested regime ISVD reduces storage by roughly two orders of magnitude relative to the direct method and by about a factor of three relative to the fast method, while retaining comparable accuracy. For long-time simulations, this difference becomes even more important because the direct method must store all $N_t$ snapshots, whereas the fast and ISVD approaches replace the full history by $N_{\exp}$ exponential modes or by low-rank factors, respectively.}
\subsection{Example 3}
In this section, we now extend the problem to three dimensions. We solve the 3D sub-diffusion equation over the unit cube domain $\Omega = (0,1)^3$. The exact solution is given by
\[
u(\mathbf{x},t) = t^{2+\alpha}x(1-x)y(1-y)z(1-z)
\]
where $f$ can be derived analytically. 

{\color{blue}Since the direct method is already memory-prohibitive for this three-dimensional problem, we compare only the ISVD and fast evaluation methods. To address the concern that the previous version used only one very small time step, we now report two complementary three-dimensional studies: first, a spatial-convergence comparison for the larger time step $\tau=10^{-4}$ and several fractional orders; second, a parameter study for the representative case $\alpha=0.5$ that varies both $\tau$ and $(\mathtt{tol1},\mathtt{tol2},\mathtt{tol3})$ and records the corresponding storage requirements.}

\begin{table}[H]
	\centering
	\color{blue}
	\caption{Spatial convergence rates and CPU times of the ISVD and Fast evaluation methods in three dimensions for $\tau=10^{-4}$ and different fractional orders $\alpha$.}
	\label{tab:3d_spatial_grouped_dof_new}
	\setlength{\tabcolsep}{4.2pt}
	\renewcommand{\arraystretch}{1.15}
	\small
	\begin{adjustbox}{max width=\textwidth}
		\begin{tabular}{c ccc ccc}
			\toprule
			& \multicolumn{3}{c}{ISVD} & \multicolumn{3}{c}{Fast Evaluation} \\
			\cmidrule(lr){2-4}\cmidrule(lr){5-7}
			$h$ & Error & Rate & CPU(s) & Error & Rate & CPU(s) \\
			\midrule
			
			\multicolumn{7}{c}{$\alpha = 0.2$} \\
			\midrule
			$1/4$  & $1.50 \times 10^{-3}$ & --   & 104.08  & $1.50 \times 10^{-3}$ & --   & 7.09 \\
			$1/8$  & $4.21 \times 10^{-4}$ & 1.83 & 103.41  & $4.21 \times 10^{-4}$ & 1.83 & 42.30 \\
			$1/16$ & $1.08 \times 10^{-4}$ & 1.96 & 463.00  & $1.08 \times 10^{-4}$ & 1.96 & 597.57 \\
			$1/32$ & $2.73 \times 10^{-5}$ & 1.99 & 3495.83 & $2.73 \times 10^{-5}$ & 1.99 & 7497.80 \\
			
			\midrule
			\multicolumn{7}{c}{$\alpha = 0.5$} \\
			\midrule
			$1/4$  & $1.51 \times 10^{-3}$ & --   & 65.19   & $1.51 \times 10^{-3}$ & --   & 6.79 \\
			$1/8$  & $4.22 \times 10^{-4}$ & 1.84 & 96.85   & $4.21 \times 10^{-4}$ & 1.84 & 41.99 \\
			$1/16$ & $1.09 \times 10^{-4}$ & 1.96 & 500.31  & $1.09 \times 10^{-4}$ & 1.96 & 550.95 \\
			$1/32$ & $2.73 \times 10^{-5}$ & 1.99 & 3515.21 & $2.73 \times 10^{-5}$ & 1.99 & 7001.66 \\
			
			\midrule
			\multicolumn{7}{c}{$\alpha = 0.8$} \\
			\midrule
			$1/4$  & $1.51 \times 10^{-3}$ & --   & 78.56   & $1.51 \times 10^{-3}$ & --   & 6.38 \\
			$1/8$  & $4.23 \times 10^{-4}$ & 1.84 & 135.26  & $4.23 \times 10^{-4}$ & 1.84 & 42.17 \\
			$1/16$ & $1.09 \times 10^{-4}$ & 1.96 & 569.38  & $1.09 \times 10^{-4}$ & 1.96 & 571.16 \\
			$1/32$ & $2.74 \times 10^{-5}$ & 1.99 & 3795.66 & $2.74 \times 10^{-5}$ & 1.99 &  \multicolumn{1}{c}{7508.20} \\
			
			\bottomrule
		\end{tabular}
	\end{adjustbox}
	
	\vspace{1mm}
	\parbox{0.97\textwidth}{\footnotesize
		Note: All errors are measured in the $L^2$ norm. Convergence rates are computed between two consecutive mesh refinements.
	}
	\end{table}

{\color{blue}Table~\ref{tab:3d_spatial_grouped_dof_new} shows that the ISVD and fast evaluation methods produce essentially identical $L^2$ errors and nearly second-order spatial convergence for all tested fractional orders. The fast evaluation CPU times increase regularly as the mesh is refined because, for fixed $\tau$, its per-step history cost is governed mainly by the spatial size and the fixed number of exponential terms. The ISVD timings are not strictly monotone on the two coarsest meshes because there the retained rank is very small and the total running time is influenced noticeably by fixed initialization and update overhead. Once the mesh is refined to $h=1/16$ and $h=1/32$, the expected scaling becomes clear; on these finer meshes ISVD is comparable to or faster than the fast evaluation method in all tested cases.}

\begin{table}[htbp]
	\centering
		\color{blue}
	\caption{CPU time and storage comparison of the ISVD and Fast methods on two 3D meshes for different time step sizes $\tau$ and ISVD tolerances with $\alpha=0.5$.}
	\label{tab:isvd_fast_tol_tau_combined}
	\setlength{\tabcolsep}{5pt}
	\renewcommand{\arraystretch}{1.15}
	\small
	\begin{adjustbox}{max width=\textwidth}
		\begin{tabular}{c c c c c c}
			\toprule
			$h$ & $\tau$ & Method & $(\mathrm{tol}_1,\mathrm{tol}_2,\mathrm{tol}_3)$ & CPU(s) & Storage (MB) \\
			\midrule
			
			\multicolumn{6}{c}{DOF$=4913$} \\
			\midrule
			\multirow{8}{*}{$1/16$}
			& \multirow{4}{*}{$2.0\times 10^{-4}$}
			& Fast & -- & 26.85 & 2.063 \\
			& & ISVD & $(10^{-8},10^{-10},10^{-6})$   & 21.69 & 0.303 \\
			& & ISVD & $(10^{-10},10^{-12},10^{-8})$  & 21.74 & 0.303 \\
			& & ISVD & $(10^{-12},10^{-14},10^{-10})$ & 21.84 & 0.454 \\
			\cmidrule(lr){2-6}
			& \multirow{4}{*}{$1.0\times 10^{-4}$}
			& Fast & -- & 56.56 & 2.138 \\
			& & ISVD & $(10^{-8},10^{-10},10^{-6})$   & 43.49 & 0.455 \\
			& & ISVD & $(10^{-10},10^{-12},10^{-8})$  & 43.78 & 0.455 \\
			& & ISVD & $(10^{-12},10^{-14},10^{-10})$ & 44.00 & 0.683 \\
			\midrule
			
			\multicolumn{6}{c}{DOF$=35937$} \\
			\midrule
			\multirow{8}{*}{$1/32$}
			& \multirow{4}{*}{$2.0\times 10^{-4}$}
			& Fast & -- & 399.94 & 15.081 \\
			& & ISVD & $(10^{-8},10^{-10},10^{-6})$   & 184.49 & 1.249 \\
			& & ISVD & $(10^{-10},10^{-12},10^{-8})$  & 193.81 & 1.249 \\
			& & ISVD & $(10^{-12},10^{-14},10^{-10})$ & 194.98 & 1.874 \\
			\cmidrule(lr){2-6}
			& \multirow{4}{*}{$1.0\times 10^{-4}$}
			& Fast & -- & 803.80 & 15.630 \\
			& & ISVD & $(10^{-8},10^{-10},10^{-6})$   & 367.46 & 1.402 \\
			& & ISVD & $(10^{-10},10^{-12},10^{-8})$  & 381.29 & 1.402 \\
			& & ISVD & $(10^{-12},10^{-14},10^{-10})$ & 383.65 & 2.103 \\
			\bottomrule
		\end{tabular}
	\end{adjustbox}
	
	\vspace{1mm}
	\parbox{0.98\textwidth}{\footnotesize
		Note: The Fast method does not involve ISVD truncation tolerances, so these entries are marked by ``--''.
	}
\end{table}

{\color{blue}To examine the influence of larger time steps and the tolerance choice more directly, Table~\ref{tab:isvd_fast_tol_tau_combined} reports additional three-dimensional results for $\tau=2\times 10^{-4}$ and $\tau=10^{-4}$ on the two finer meshes $h=1/16$ and $h=1/32$ with $\alpha=0.5$. Increasing $\tau$ reduces the total CPU time for both methods because the number of time levels is halved, but the table shows that ISVD remains efficient across all tested tolerance sets. The storage comparison is explicit: when $h=1/32$ and $\tau=10^{-4}$, the fast evaluation method requires $15.630$~MB, while ISVD uses only $1.402$--$2.103$~MB; when $\tau=2\times 10^{-4}$, the corresponding values are $15.081$~MB for the fast method and $1.249$--$1.874$~MB for ISVD. Taken together with Table~\ref{tab:3d_spatial_grouped_dof_new}, these results show that the revised three-dimensional study now documents the influence of both $\tau$ and the ISVD tolerances much more explicitly, while also demonstrating clear storage savings and competitive or smaller CPU times on the finer meshes.}

\section{Conclusion} 
		{\color{blue} In this paper, we propose an efficient and memory-saving algorithm for solving time-fractional PDEs by using  incremental SVD algorithm. The numerical results show that the method preserves the expected accuracy, dramatically reduces memory usage relative to the direct method, and remains competitive with a representative fast evaluation method in the tested parameter regimes. Although we only test the fractional heat equation, the same framework can potentially be extended to other time-fractional PDEs and to models with substantial history-storage demands, such as Cole--Cole models for Maxwell's equations  \cite{Li_SISC2011,FAN2025116566}. The present analysis is normwise: when the truncation tolerances are sufficiently small, the compressed-history solution remains close to the standard fully discrete solution in the $L^2$ sense, but we do not prove preservation of structural properties such as positivity, monotonicity, or a discrete maximum principle. Developing structure-preserving variants and extending the framework to other time-dependent PDEs with memory effects are important directions for future work  \cite{Zhang2011,Santiago2017,Sadiq2026,Jesus2020}.}
\section*{CRediT Authorship Contribution Statement} 
\textbf{Jichun Li}: Conceptualization, Supervision, Writing. \textbf{Yangpeng Zhang}: Writing, Software, Methodology, Investigation. \textbf{Yangwen Zhang}: Writing, Conceptualization, Methodology, Investigation, Validation.

\section*{Funding}
Yangwen Zhang is supported by the National Science Foundation DMS-2111315.

\section*{Acknowledgments}
The authors would like to express sincere gratitude to Prof. Jiwei Zhang for generously providing the source code for the fast evaluation method.
\bibliographystyle{siamplain}
\bibliography{references_AMS_cleaned}

\begin{thebibliography}{10}

\bibitem{Santiago2017}
{\sc S.~Badia and J.~Bonilla}, {\em Monotonicity-preserving finite element
  schemes based on differentiable nonlinear stabilization}, Comput. Methods
  Appl. Mech. Engrg., 313 (2017), pp.~133--158,
  \url{https://doi.org/10.1016/j.cma.2016.09.035}.

\bibitem{BanjaiLopez2019}
{\sc L.~Banjai and M.~L\'{o}pez-Fern\'{a}ndez}, {\em Efficient high order
  algorithms for fractional integrals and fractional differential equations},
  Numer. Math., 141 (2019), pp.~289--317,
  \url{https://doi.org/10.1007/s00211-018-1004-0}.

\bibitem{Benson2000}
{\sc D.~A. Benson, S.~W. Wheatcraft, and M.~M. Meerschaert}, {\em Application
  of a fractional advection-dispersion equation}, Water Resour. Res., 36
  (2000), pp.~1403--1412.

\bibitem{Jesus2020}
{\sc J.~Bonilla and S.~Badia}, {\em Monotonicity-preserving finite element
  schemes with adaptive mesh refinement for hyperbolic problems}, J. Comput.
  Phys., 416 (2020), pp.~109522, 23,
  \url{https://doi.org/10.1016/j.jcp.2020.109522}.

\bibitem{Brand2006}
{\sc M.~Brand}, {\em Fast low-rank modifications of the thin singular value
  decomposition}, Linear Algebra Appl., 415 (2006), pp.~20--30,
  \url{https://doi.org/10.1016/j.laa.2005.07.021}.

\bibitem{CheniFEM}
{\sc L.~Chen}, {\em {iFEM}: an integrated finite element methods package in
  {MATLAB}}, tech. report, University of California, Irvine, 2009.

\bibitem{MR4982783}
{\sc E.~Di~Costanzo, N.~K\"{u}hl, J.-C. Marongiu, and T.~Rung}, {\em
  Incremental model order reduction of smoothed-particle hydrodynamic
  simulations}, Internat. J. Numer. Methods Fluids, 97 (2025), pp.~1571--1594,
  \url{https://doi.org/10.1002/fld.70012}.

\bibitem{FAN2025116566}
{\sc E.~Fan, J.~Li, Y.~Liu, and Y.~Zhang}, {\em Isogeometric analysis for
  time-dependent {M}axwell's equations in complex media}, J. Comput. Appl.
  Math., 465 (2025), pp.~Paper No. 116566, 17,
  \url{https://doi.org/10.1016/j.cam.2025.116566}.

\bibitem{Areed2020}
{\sc H.~Fareed and J.~R. Singler}, {\em Error analysis of an incremental proper
  orthogonal decomposition algorithm for {PDE} simulation data}, J. Comput.
  Appl. Math., 368 (2020), pp.~112525, 14,
  \url{https://doi.org/10.1016/j.cam.2019.112525}.

\bibitem{MR3775096}
{\sc H.~Fareed, J.~R. Singler, Y.~Zhang, and J.~Shen}, {\em Incremental proper
  orthogonal decomposition for {PDE} simulation data}, Comput. Math. Appl., 75
  (2018), pp.~1942--1960, \url{https://doi.org/10.1016/j.camwa.2017.09.012}.

\bibitem{Garrappa2018}
{\sc R.~Garrappa}, {\em Numerical solution of fractional differential
  equations: a survey and a software tutorial}, Mathematics, 6 (2018), p.~16.

\bibitem{Giraud2005}
{\sc L.~Giraud, J.~Langou, M.~Rozlo\v{z}n\'{\i}k, and J.~van~den Eshof}, {\em
  Rounding error analysis of the classical {G}ram-{S}chmidt orthogonalization
  process}, Numer. Math., 101 (2005), pp.~87--100,
  \url{https://doi.org/10.1007/s00211-005-0615-4}.

\bibitem{Sadiq2026}
{\sc S.~Hamidi, M.~{El Ossmani}, and A.~Taakili}, {\em A nonlinear finite
  volume scheme preserving positivity for the numerical simulation of
  groundwater contaminant transport in porous media}, J. Comput. Sci., 93
  (2026), p.~102739, \url{https://doi.org/10.1016/j.jocs.2025.102739}.

\bibitem{Zhuang2019variableCoeff}
{\sc L.~He and J.~Lv}, {\em Efficient finite element numerical solution of the
  variable coefficient fractional subdiffusion equation}, Adv. Difference Equ.,
   (2019), pp.~Paper No. 116, 17,
  \url{https://doi.org/10.1186/s13662-019-2048-x}.

\bibitem{Jiang2015}
{\sc S.~Jiang, L.~Greengard, and S.~Wang}, {\em Efficient sum-of-exponentials
  approximations for the heat kernel and their applications}, Adv. Comput.
  Math., 41 (2015), pp.~529--551,
  \url{https://doi.org/10.1007/s10444-014-9372-1}.

\bibitem{Zhang_fast_evaluation2017}
{\sc S.~Jiang, J.~Zhang, Q.~Zhang, and Z.~Zhang}, {\em Fast evaluation of the
  {C}aputo fractional derivative and its applications to fractional diffusion
  equations}, Commun. Comput. Phys., 21 (2017), pp.~650--678,
  \url{https://doi.org/10.4208/cicp.OA-2016-0136}.

\bibitem{Jin2014}
{\sc B.~Jin, R.~Lazarov, J.~Pasciak, and Z.~Zhou}, {\em Error analysis of
  semidiscrete finite element methods for inhomogeneous time-fractional
  diffusion}, IMA J. Numer. Anal., 35 (2015), pp.~561--582,
  \url{https://doi.org/10.1093/imanum/dru018}.

\bibitem{Jin2013}
{\sc B.~Jin, R.~Lazarov, and Z.~Zhou}, {\em Error estimates for a semidiscrete
  finite element method for fractional order parabolic equations}, SIAM J.
  Numer. Anal., 51 (2013), pp.~445--466,
  \url{https://doi.org/10.1137/120873984}.

\bibitem{Jin2015AN}
{\sc B.~Jin, R.~Lazarov, and Z.~Zhou}, {\em An analysis of the {L}1 scheme for
  the subdiffusion equation with nonsmooth data}, IMA J. Numer. Anal., 36
  (2016), pp.~197--221, \url{https://doi.org/10.1093/imanum/dru063}.

\bibitem{JinLiZhou2016CN}
{\sc B.~Jin, B.~Li, and Z.~Zhou}, {\em An analysis of the {C}rank-{N}icolson
  method for subdiffusion}, IMA J. Numer. Anal., 38 (2018), pp.~518--541,
  \url{https://doi.org/10.1093/imanum/drx019}.

\bibitem{Kilbas2006}
{\sc A.~A. Kilbas, H.~M. Srivastava, and J.~J. Trujillo}, {\em Theory and
  applications of fractional differential equations}, vol.~204 of North-Holland
  Mathematics Studies, Elsevier Science B.V., Amsterdam, 2006.

\bibitem{LarsonBengzon2013}
{\sc M.~G. Larson and F.~Bengzon}, {\em The finite element method: theory,
  implementation, and applications}, vol.~10 of Texts in Computational Science
  and Engineering, Springer, Heidelberg, 2013,
  \url{https://doi.org/10.1007/978-3-642-33287-6}.

\bibitem{Le2016}
{\sc K.~N. Le, W.~McLean, and K.~Mustapha}, {\em Numerical solution of the
  time-fractional {F}okker-{P}lanck equation with general forcing}, SIAM J.
  Numer. Anal., 54 (2016), pp.~1763--1784,
  \url{https://doi.org/10.1137/15M1031734}.

\bibitem{LiYangZhou2023splitFEM}
{\sc B.~Li, Z.~Yang, and Z.~Zhou}, {\em High-order splitting finite element
  methods for the subdiffusion equation with limited smoothing property}, Math.
  Comp., 93 (2024), pp.~2557--2586, \url{https://doi.org/10.1090/mcom/3944}.

\bibitem{Li_SISC2011}
{\sc J.~Li, Y.~Huang, and Y.~Lin}, {\em Developing finite element methods for
  {M}axwell's equations in a {C}ole-{C}ole dispersive medium}, SIAM J. Sci.
  Comput., 33 (2011), pp.~3153--3174, \url{https://doi.org/10.1137/110827624}.

\bibitem{Li2010fasttime}
{\sc J.-R. Li}, {\em A fast time stepping method for evaluating fractional
  integrals}, SIAM J. Sci. Comput., 31 (2009/10), pp.~4696--4714,
  \url{https://doi.org/10.1137/080736533}.

\bibitem{li2024incremental}
{\sc X.~Li, J.~R. Singler, and X.~He}, {\em Incremental data compression for
  {PDE}-constrained optimization with a data assimilation application}, arXiv
  preprint arXiv:2404.09323,  (2024).

\bibitem{LIN20071533}
{\sc Y.~Lin and C.~Xu}, {\em Finite difference/spectral approximations for the
  time-fractional diffusion equation}, J. Comput. Phys., 225 (2007),
  pp.~1533--1552, \url{https://doi.org/10.1016/j.jcp.2007.02.001}.

\bibitem{McLean2012interval}
{\sc W.~McLean}, {\em Fast summation by interval clustering for an evolution
  equation with memory}, SIAM J. Sci. Comput., 34 (2012), pp.~A3039--A3056,
  \url{https://doi.org/10.1137/120870505}.

\bibitem{METZLER20001}
{\sc R.~Metzler and J.~Klafter}, {\em The random walk's guide to anomalous
  diffusion: a fractional dynamics approach}, Phys. Rep., 339 (2000), p.~77,
  \url{https://doi.org/10.1016/S0370-1573(00)00070-3}.

\bibitem{Mustapha2018}
{\sc K.~Mustapha}, {\em F{EM} for time-fractional diffusion equations, novel
  optimal error analyses}, Math. Comp., 87 (2018), pp.~2259--2272,
  \url{https://doi.org/10.1090/mcom/3304}.

\bibitem{Ren2017superconvergence}
{\sc J.~Ren, X.~Long, S.~Mao, and J.~Zhang}, {\em Superconvergence of finite
  element approximations for the fractional diffusion-wave equation}, J. Sci.
  Comput., 72 (2017), pp.~917--935,
  \url{https://doi.org/10.1007/s10915-017-0385-z}.

\bibitem{Ren2018fast}
{\sc J.~Ren, S.~Mao, and J.~Zhang}, {\em Fast evaluation and high accuracy
  finite element approximation for the time fractional subdiffusion equation},
  Numer. Methods Partial Differential Equations, 34 (2018), pp.~705--730,
  \url{https://doi.org/10.1002/num.22226}.

\bibitem{SakamotoYamamoto2011}
{\sc K.~Sakamoto and M.~Yamamoto}, {\em Initial value/boundary value problems
  for fractional diffusion-wave equations and applications to some inverse
  problems}, J. Math. Anal. Appl., 382 (2011), pp.~426--447,
  \url{https://doi.org/10.1016/j.jmaa.2011.04.058}.

\bibitem{Schadle2006}
{\sc A.~Sch\"{a}dle, M.~L\'{o}pez-Fern\'{a}ndez, and C.~Lubich}, {\em Fast and
  oblivious convolution quadrature}, SIAM J. Sci. Comput., 28 (2006),
  pp.~421--438, \url{https://doi.org/10.1137/050623139}.

\bibitem{Scher1975}
{\sc H.~Scher and E.~W. Montroll}, {\em Anomalous transit-time dispersion in
  amorphous solids}, Phys. Rev. B, 12 (1975), pp.~2455--2477.

\bibitem{Zeng2013}
{\sc F.~Zeng, C.~Li, F.~Liu, and I.~Turner}, {\em The use of finite
  difference/element approaches for solving the time-fractional subdiffusion
  equation}, SIAM J. Sci. Comput., 35 (2013), pp.~A2976--A3000,
  \url{https://doi.org/10.1137/130910865}.

\bibitem{Zeng2015second}
{\sc F.~Zeng, C.~Li, F.~Liu, and I.~Turner}, {\em Numerical algorithms for
  time-fractional subdiffusion equation with second-order accuracy}, SIAM J.
  Sci. Comput., 37 (2015), pp.~A55--A78,
  \url{https://doi.org/10.1137/14096390X}.

\bibitem{Zeng2016fasthighdim}
{\sc F.~Zeng, Z.~Zhang, and G.~E. Karniadakis}, {\em Fast difference schemes
  for solving high-dimensional time-fractional subdiffusion equations}, J.
  Comput. Phys., 307 (2016), pp.~15--33,
  \url{https://doi.org/10.1016/j.jcp.2015.11.058}.

\bibitem{Zhang2011}
{\sc X.~Zhang and C.-W. Shu}, {\em Maximum-principle-satisfying and
  positivity-preserving high-order schemes for conservation laws: survey and
  new developments}, Proc. R. Soc. Lond. Ser. A Math. Phys. Eng. Sci., 467
  (2011), pp.~2752--2776, \url{https://doi.org/10.1098/rspa.2011.0153}.

\bibitem{Zhang2022isvd}
{\sc Y.~Zhang}, {\em An answer to an open question in the incremental {SVD}},
  2022.
\newblock arXiv:2204.05398.

\end{thebibliography}

\end{document}